\documentclass{amsart}
\usepackage{graphicx}
\usepackage{amscd}
\usepackage{amsmath}
\usepackage{amsfonts}
\usepackage{amssymb}
\newtheorem{theorem}{Theorem}[section]

\newtheorem{corollary}[theorem]{Corollary}

\newtheorem{lemma}[theorem]{Lemma}

\newtheorem{proposition}[theorem]{Proposition}
\newtheorem{remark}[theorem]{Remark}

\numberwithin{equation}{section}

\begin{document}
\title[Maps with prescribed Thom-Boardman singularities]{A homotopy principle for maps with prescribed Thom-Boardman singularities}
\author{YOSHIFUMI ANDO}
\address{Department of Mathematical Sciences, Faculty of Science, Yamaguchi University,
Yamaguchi 753-8512, Japan}
\email{andoy@yamaguchi-u.ac.jp}
\date{}
\dedicatory{}

\begin{abstract}
Let $N$ and $P$ be smooth manifolds of dimensions $n$ and $p$ ($n\geq p\geq2$)
respectively. Let $\Omega^{I}(N,P)$ denote an open subspace of $J^{\infty
}(N,P)$ which consists of all Boardman submanifolds $\Sigma^{J}(N,P)$ of
symbols $J$ with $J\leq I$. An $\Omega^{I}$-regular map $f:N\rightarrow P$
refers to a smooth map having only singularities in $\Omega^{I}(N,P)$ and
satisfying transversality condition. We will prove what is called the homotopy
principle for $\Omega^{I}$-regular maps in the existence level. Namely, a
continuous section $s$ of $\Omega^{I}(N,P)$ over $N$ has an $\Omega^{I}%
$-regular map $f$ such that $s$ and $j^{\infty}f$ are homotopic as sections.
\end{abstract}\maketitle

\section*{Introduction}

Let $N$ and $P$ be smooth ($C^{\infty}$) manifolds of dimensions $n$ and $p$
respectively with $n\geq p\geq2$. In [B] there have been defined what are
called the Boardman manifolds $\Sigma^{I}(N,P)$ in $J^{\infty}(N,P)$ for the
symbol $I=(i_{1},i_{2},\cdots,i_{r})$, where $i_{1},i_{2},\cdots,i_{r}$ are a
finite number of integers with $i_{1}\geq i_{2}\geq\cdots\geq i_{r}\geq0$. We
say that a smooth map germ $f:(N,x)\rightarrow(P,y)$ has $x$ as a
Thom-Boardman singularity of the symbol $I$ if and only if $j_{x}^{\infty}%
f\in\Sigma^{I}(N,P)$. Let $\Omega^{I}(N,P)$ denote an open subset of
$J^{\infty}(N,P)$ which consists of all Boardman manifolds $\Sigma^{J}(N,P)$
with symbols $J$ of length $r$ satisfying $J\leq I$ in the lexicographic
order. It is known that $\Omega^{I}(N,P)$ is an open subbundle of $J^{\infty
}(N,P)$ with the projection $\pi_{N}^{\infty}\times\pi_{P}^{\infty}$, whose
fiber is denoted by $\Omega^{I}(n,p)$. A smooth map $f:N\rightarrow P$ is
called an $\Omega^{I}$-\textit{regular map} if and only if (i) $j^{\infty
}f(N)\subset\Omega^{I}(N,P)$ and (ii) $j^{\infty}f$ is transverse to all
$\Sigma^{J}(N,P)$.

We will study a homotopy theoretic condition for a given continuous map to be
homotopic to an $\Omega^{I}$-regular map. Let $C_{\Omega^{I}}^{\infty}(N,P)$
denote the space consisting of all $\Omega^{I}$-regular maps equipped with the
$C^{\infty}$-topology. Let $\Gamma_{\Omega^{I}}(N,P)$ denote the space
consisting of all continuous sections of the fiber bundle $\pi_{N}^{\infty
}|\Omega^{I}(N,P):\Omega^{I}(N,P)\rightarrow N$ equipped with the compact-open
topology. Then there exists a continuous map
\[
j_{\Omega^{I}}:C_{\Omega^{I}}^{\infty}(N,P)\rightarrow\Gamma_{\Omega^{I}%
}(N,P)
\]
defined by $j_{\Omega^{I}}(f)=j^{\infty}f$. It follows from the well-known
theorem due to Gromov[G1] that if $N$ is a connected open manifold, then
$j_{\Omega^{I}}$ is a weak homotopy equivalence. This property is called the
homotopy principle (the terminology used in [G2]). If $N$ is a closed
manifold, then it becomes a hard problem for us to prove the homotopy
principle. As the primary investigation preceding [G1], we must refer to the
Smale-Hirsch Immersion Theorem ([H]), $k$-mersion Theorem due to [F] and the
Phillips Submersion Theorem for open manifolds ([P]). In [E1] and [E2],
\`{E}lia\v{s}berg has proved the well-known homotopy principle in the $1$-jet
level for $\Omega^{n-p+1,0}$-regular maps, say fold-maps. As for the
Thom-Boardman singularities, du Plessis[duP] has proved that if $i_{r}%
>n-p-d^{I}$, where $d^{I}$ is the sum of $\alpha_{1},\cdots,\alpha_{r-1}$ with
$\alpha_{\ell}$ being $1$ or $0$ depending on $i_{\ell}-i_{\ell+1}>1$ or
otherwise, then $j_{\Omega^{I}}$ is a weak homotopy equivalence.

In this paper we prove the following homotopy principle in the existence level
for closed manifolds.

\begin{theorem}
Let $n\geq p\geq2$. Let $N$ and $P$ be connected manifolds of dimensions $n$
and $p$ respectively with $\partial N=\emptyset$. Assume that $\Omega
^{I}(N,P)$ contains $\Sigma^{n-p+1,0}(N,P)$ at least. Let $C$ be a closed
subset of $N$. Let $s$ be a section of $\Gamma_{\Omega^{I}}(N,P)$ which has an
$\Omega^{I}$-regular map $g$ defined on a neighborhood of $C$ into $P$, where
$j^{\infty}g=s$.

Then there exists an $\Omega^{I}$-regular map $f:N\rightarrow P$ such that
$j^{\infty}f$ is homotopic to $s$ relative to a neighborhood of $C$ by a
homotopy $s_{\lambda}$ in $\Gamma_{\Omega^{I}}(N,P)$ with $s_{0}=s$ and
$s_{1}=j^{\infty}f$.
\end{theorem}

In [A1] we have given Theorem 0.1 for the symbol $I=(n-p+1,\overbrace
{1,\cdots,1}^{r-1},0)$ with a partially sketchy proof using the results in
[E1] and [E2]. The singularities of this symbol $I$ are often called $A_{r}%
$-singularities or Morin singularities. The detailed proof are given in [An4,
Theorem 4.1] and [An6, Theorem 0.5] for the symbol $I=(n-p+1,0)$. We will use
these two theorems in the proof of Theorem 0.1 in this paper.

Recently it turns out that this kind of the homotopy principle has many
applications. Theorem 0.1 is very important even for fold-maps in proving the
relations between fold-maps, surgery theory and stable homotopy groups ([An4,
Theorem 1] and [An5, Theorems 0.2 and 0.3]). The homotopy type of
$\Omega^{n-p+1,0}$ determined in [An3] and [An5] has played an important role.
We can now readily deduce the famous theorem about the elimination of cusps in
[L1] and [E1] (see also [T]) from these theorems.

The homotopy principle in the existence level for maps and singular foliations
having only what are called $A$, $D$ and $E$ singularities are proved in [An2]
and [An7].

In [Sady] Sadykov has applied [An1, Theorem 1] to the elimination of higher
$A_{r}$ singularities ($r\geq3$) for Morin maps when $n-p$ is odd. This result
is a strengthened version of the Chess conjecture proposed in [C].

As an application of Theorem 0.1 we prove the following theorem. We note that
the simplest case is also a little stronger form of the Chess conjecture.

\begin{theorem}
Let $n\geq p\geq2,$ and $N$ and $P$ be connected manifolds of dimensions $n$
and $p$ respectively. Let $I=(n-p+1,i_{2},\cdots,i_{r-1},1,1)$ and
$J=(n-p+1,i_{2},\cdots,i_{r-1},1,0)$ such that $n-p+1-i_{2}$ and $r$
$(r\geq3)$ are odd integers. Then if $f:N\rightarrow P$ is an $\Omega^{I}%
$-regular map, then $f$ is homotopic to an $\Omega^{J}$-regular map
$g:N\rightarrow P$ such that $j^{\infty}f$ and $j^{\infty}g$ are homotopic in
$\Gamma_{\Omega^{I}}(N,P)$.
\end{theorem}

In Section 1 we explain notations which are used in this paper. In Section
2\ we review the definitions and the fundamental properties of the Boardman
manifolds, from which we deduce several further results about higher intrinsic
derivatives in Section 3. In Section 4 we reduce the proof of Theorem 0.1 to
the proof of Theorem 4.1 by the induction, and prepare a certain rotation of
the tangent spaces defined around the singularities of given symbol in $N$ to
deform the section $s$. In Section 5 we prepare several lemmas which are used
in the deformation of the section $s$ in the proof of Theorem 4.1. We prove
Theorem 4.1 in Section 6 and prove Theorem 0.2 in Section 7.

\section{Notations}

Throughout the paper all manifolds are Hausdorff, paracompact and smooth of
class $C^{\infty}$. Maps are basically continuous, but may be smooth (of class
$C^{\infty}$) if necessary. Given a fiber bundle $\pi:E\rightarrow X$ and a
subset $C$ in $X,$ we denote $\pi^{-1}(C)$ by $E|_{C}.$ Let $\pi^{\prime
}:F\rightarrow Y$ be another fiber bundle. A map $\tilde{b}:E\rightarrow F$ is
called a fiber map over a map $b:X\rightarrow Y$ if $\pi^{\prime}\circ
\tilde{b}=b\circ\pi$ holds. The restriction $\tilde{b}|(E|_{C}):E|_{C}%
\rightarrow F$ (or $F|_{b(C)}$) is denoted by $\tilde{b}|_{C}$. In particular,
for a point $x\in X,$ $E|_{x}$ and $\tilde{b}|_{x}$ are simply denoted by
$E_{x}$ and $\tilde{b}_{x}:E_{x}\rightarrow F_{b(x)}$ respectively. We denote,
by $b^{F}$, the induced fiber map $b^{\ast}(F)\rightarrow F$ covering $b$. For
a map $j:W\rightarrow X$, let $j^{\ast}(\tilde{b}):j^{\ast}E\rightarrow(b\circ
j)^{\ast}F$ over $W$ be the fiber map canonically induced from $b$ and $j$. A
fiberwise homomorphism $E\rightarrow F$ is simply called a homomorphism. For a
vector bundle $E$ with a metric and a positive function $\delta$ on $X$, let
$D_{\delta}(E)$ be the associated disk bundle of $E$ with radius $\delta$. If
there is a canonical isomorphism between two vector bundles $E$ and $F$ over
$X=Y,$ then we write $E\cong F$.

When $E$ and $F$ are smooth vector bundles over $X=Y$, Hom$(E,F)$ denotes the
smooth vector bundle over $X$ with fiber Hom$(E_{x},F_{x})$, $x\in X$ which
consists of all homomorphisms $E_{x}\rightarrow F_{x}$.

Let $J^{k}(N,P)$ denote the $k$-jet space of manifolds $N$ and $P$. Let
$\pi_{N}^{k}$ and $\pi_{P}^{k}$ be the projections mapping a jet to its source
and target respectively. The map $\pi_{N}^{k}\times\pi_{P}^{k}:J^{k}%
(N,P)\rightarrow N\times P$ induces a structure of a fiber bundle with
structure group $L^{k}(p)\times L^{k}(n)$, where $L^{k}(m)$ denotes the group
of all $k$-jets of local diffeomorphisms of $(\mathbf{R}^{m},0)$. The fiber
$(\pi_{N}^{k}\times\pi_{P}^{k})^{-1}(x,y)$ is denoted by $J_{x,y}^{k}(N,P)$.

Let $\pi_{N}$ and $\pi_{P}$ be the projections of $N\times P$ onto $N$ and $P$
respectively.\ We set
\begin{equation}
J^{k}(TN,TP)=\bigoplus_{i=1}^{k}\text{\textrm{Hom}}(S^{i}(\pi_{N}^{\ast
}(TN)),\pi_{P}^{\ast}(TP))
\end{equation}
over $N\times P$. Here, for a vector bundle $E$ over $X$, let $S^{i}(E)$ be
the vector bundle $\cup_{x\in X}S^{i}(E_{x})$ over $X$, where $S^{i}(E_{x})$
denotes the $i$-fold symmetric product of $E_{x}$. If we provide $N$ and $P$
with Riemannian metrics, then the Levi-Civita connections induce the
exponential maps $\exp_{N,x}:T_{x}N\rightarrow N$ and $\exp_{P,y}%
:T_{y}P\rightarrow P$. In dealing with the exponential maps we always consider
the convex neighborhoods ([K-N]). We define the smooth bundle map
\begin{equation}
J^{k}(N,P)\mathbf{\rightarrow}J^{k}(TN,TP)\text{ \ \ \ over }N\times P
\end{equation}
by sending $z=j_{x}^{k}f\in J_{x,y}^{k}(N,P)$ to the $k$-jet of $(\exp
_{P,y})^{-1}\circ f\circ\exp_{N,x}$ at $\mathbf{0}\in T_{x}N$, which is
regarded as an element of $J^{k}(T_{x}N,T_{y}P)(=J_{x,y}^{k}(TN,TP))$ (see
[K-N, Proposition 8.1] for the smoothness of exponential maps). More strictly,
(1.2) gives a smooth equivalence of the fiber bundles under the structure
group $L^{k}(p)\times L^{k}(n)$. Namely, it gives a smooth reduction of the
structure group $L^{k}(p)\times L^{k}(n)$ of $J^{k}(N,P)$ to $O(p)\times
O(n)$, which is the structure group of $J^{k}(TN,TP)$.

Recall that $S^{i}(E)\ $has the inclusion $S^{i}(E)\rightarrow\otimes^{i}E$
and the canonical projection $\otimes^{i}E\rightarrow S^{i}(E)$ (see [B,
Section 4] and [Mats, Ch. III, Section 2]). Let $E_{j}$ be subbundles of $E$
$(j=1,\cdots,i)$. We define $E_{1}\bigcirc\cdots\bigcirc E_{i}=\bigcirc
_{j=1}^{i}E_{j}$ to be the image of $E_{1}\otimes\cdots\otimes E_{i}%
=\otimes_{j=1}^{i}E_{j}\rightarrow\otimes^{i}E\rightarrow S^{i}(E)$. When
$E_{j+1}=\cdots=E_{j+\ell}$, we often write $E_{1}\bigcirc\cdots\bigcirc
E_{j}\bigcirc^{\ell}E_{j+1}\bigcirc E_{j+\ell+1}\bigcirc\cdots\bigcirc E_{i}$
in place of $\bigcirc_{j=1}^{i}E_{j}.$

\section{Boardman manifolds}

We review well-known results about Boardman manifolds in $J^{\infty}%
(N,P)$\ ([B], [L2] and [Math2]). Let $I=(i_{1},\cdots,i_{r})$ be a Boardman
symbol with $i_{1}\geq\cdots\geq i_{r}\geq0$. For $k\leq r$, set $I_{k}%
=(i_{1},i_{2},\cdots,i_{k})$ and $(I_{k},0)=(i_{1},i_{2},\cdots,i_{k},0)$. In
the infinite jet space $J^{\infty}(N,P)$, there have been defined a sequence
of the submanifolds $\Sigma^{I_{1}}(N,P)\supseteq\cdots\supseteq\Sigma^{I_{r}%
}(N,P)$ with the following properties. In this paper we often write
$\Sigma^{I_{r}}$ for $\Sigma^{I_{r}}(N,P)$ if there is no confusion.

Let $\mathbf{P}=(\pi_{P}^{\infty})^{\ast}(TP)$ and $\mathbf{D}$ be the total
tangent bundle defined over $J^{\infty}(N,P)$.\ We explain important
properties of the total tangent bundle $\mathbf{D}$, which are often used in
this paper. Let $f:(N,x)\rightarrow(P,y)$ be a germ and $\digamma$ be a smooth
function in the sense of [B, Definition 1.4] defined on a neighborhood of
$j_{x}^{\infty}f$. Given a vector field $v$ defined on a neighborhood of $x$
in $N$, there is a total vector field $D$ defined on a neighborhood of
$j_{x}^{\infty}f$ such that $D\digamma\circ j^{\infty}f=v(\digamma\circ
j^{\infty}f)$. It follows that $d(j^{\infty}f)(v)(\digamma)=D\digamma
(j^{\infty}f)$ for $d(j^{\infty}f):TN\rightarrow T(J^{\infty}(N,P))$ around
$x$. This implies $d(j^{\infty}f)(v)=D.$ Hence, we have $\mathbf{D\cong(}%
\pi_{N}^{\infty})^{\ast}(TN)$.

First we have the first derivative $\mathbf{d}_{1}:\mathbf{D}\rightarrow
\mathbf{P}$\ over $J^{\infty}(N,P)$. We define $\Sigma^{I_{1}}(N,P)$ to be the
submanifold of $J^{\infty}(N,P)$ which consists of all jets $z$ such that the
kernel rank of $\mathbf{d}_{1,z}$ is $i_{1}$. Since $\mathbf{d}_{1}%
|_{\Sigma^{I_{1}}(N,P)}$ is of constant rank $n-i_{1}$, we set $\mathbf{K}%
_{1}=$Ker$(\mathbf{d}_{1})$ and $\mathbf{Q}_{1}=$Cok$(\mathbf{d}_{1})$, which
are vector bundles over $\Sigma^{I_{1}}(N,P)$. Set $\mathbf{K}_{0}=\mathbf{D}%
$, $\mathbf{P}_{0}=\mathbf{P}$ and $\Sigma^{I_{0}}(N,P)=J^{\infty}(N,P)$. We
can inductively define $\Sigma^{I_{k}}(N,P)$ and the bundles $\mathbf{K}_{k}$
and $\mathbf{P}_{k}$ over $\Sigma^{I_{k}}(N,P)$ ($k\geq1$) with the properties:

(1) $\mathbf{K}_{k-1}|_{\Sigma^{I_{k}}(N,P)}\supseteq\mathbf{K}_{k}$ over
$\Sigma^{I_{k}}(N,P)$.

(2) $\mathbf{K}_{k}$ is an $i_{k}$-dimensional subbundle of $T(\Sigma
^{I_{k-1}}(N,P))|_{\Sigma^{I_{k}}(N,P)}.$

(3) There exists the $(k+1)$-th intrinsic derivative $\mathbf{d}%
_{k+1}:T(\Sigma^{I_{k-1}}(N,P))|_{\Sigma^{I_{k}}(N,P)}\linebreak
\rightarrow\mathbf{P}_{k}$ over $\Sigma^{I_{k}}(N,P),$ so that it induces the
exact sequence\ over $\Sigma^{I_{k}}(N,P):$%
\begin{equation}%
\begin{array}
[c]{l}%
\mathbf{0\rightarrow}T(\Sigma^{I_{k}}(N,P))\overset{\text{inclusion}%
}{\hookrightarrow}T(\Sigma^{I_{k-1}}(N,P))|_{\Sigma^{I_{k}}(N,P)}%
\overset{\mathbf{d}_{k+1}}{\longrightarrow}\mathbf{P}_{k}\text{ }%
\rightarrow\mathbf{0.}%
\end{array}
\end{equation}
Namely, $\mathbf{d}_{k+1}$ induces the isomorphism of the normal bundle%
\begin{equation}
\nu(I_{k}\subset I_{k-1})=(T(\Sigma^{I_{k-1}}(N,P))|_{\Sigma^{I_{k}}%
(N,P)})/T(\Sigma^{I_{k}}(N,P))
\end{equation}
of $\Sigma^{I_{k}}(N,P)$ in $\Sigma^{I_{k-1}}(N,P)$ onto $\mathbf{P}_{k}$.

(4) $\Sigma^{I_{k+1}}(N,P)$ is defined to be the submanifold of $\Sigma
^{I_{k}}(N,P)$ which consists of all jets $z$ with dim(Ker$(\mathbf{d}%
_{k+1,z}|\mathbf{K}_{k,z}))=i_{k+1}$. In particular, $\Sigma^{I_{k}}(N,P)$ is
the disjoint union $\cup_{j=0}^{i_{k}}\Sigma^{(I_{k},j)}(N,P).$

(5) Set $\mathbf{K}_{k+1}=$Ker$(\mathbf{d}_{k+1}|\mathbf{K}_{k})$ and
$\mathbf{Q}_{k+1}\mathbf{=}$Cok$(\mathbf{d}_{k+1}|\mathbf{K}_{k})$ over
$\Sigma^{I_{k+1}}(N,P)$. Then it follows that $(\mathbf{K}_{k}|_{\Sigma
^{I_{k+1}}(N,P)})\cap T(\Sigma^{I_{k}}(N,P))|_{\Sigma^{I_{k+1}}(N,P)}%
=\mathbf{K}_{k+1}$. We have the canonical projection $\mathbf{e}%
_{k}:\mathbf{P}_{k}|_{\Sigma^{I_{k+1}}(N,P)}\rightarrow\mathbf{Q}_{k+1}$.

(6) The intrinsic derivative%
\[
d(\mathbf{d}_{k+1}|\mathbf{K}_{k}):T(\Sigma^{I_{k}}(N,P))|_{\Sigma^{I_{k+1}%
}(N,P)}\rightarrow\text{Hom}(\mathbf{K}_{k+1},\mathbf{Q}_{k+1})\text{ \ \ over
}\Sigma^{I_{k+1}}(N,P)
\]
of $\mathbf{d}_{k+1}|\mathbf{K}_{k}$ is of constant rank $\dim(\Sigma^{I_{k}%
}(N,P))-\dim(\Sigma^{I_{k+1}}(N,P))$. We set $\mathbf{P}_{k+1}%
=\operatorname{Im}(d(\mathbf{d}_{k+1}|\mathbf{K}_{k}))$ and define
$\mathbf{d}_{k+2}$ to be%
\begin{equation}
\mathbf{d}_{k+2}=d(\mathbf{d}_{k+1}|\mathbf{K}_{k}):T(\Sigma^{I_{k}%
}(N,P))|_{\Sigma^{I_{k+1}}(N,P)}\rightarrow\mathbf{P}_{k+1}%
\end{equation}
as the epimorphism.

(7) There exists the bundle homomorphism of constant rank%
\begin{equation}
\mathbf{u}_{k}:\text{Hom}(\mathbf{K}_{k}\bigcirc\mathbf{K}_{k-1}\bigcirc
\cdots\bigcirc\mathbf{K}_{1},\mathbf{P})\rightarrow\text{Hom}(\mathbf{K}%
_{k},\mathbf{Q}_{k})\text{ \ \ over }\Sigma^{I_{k}}(N,P)
\end{equation}
such that the image of $\mathbf{u}_{k}$ coincides with $\mathbf{P}_{k}$. We
denote, by $\mathbf{c}_{k}$, the map $\mathbf{u}_{k}$ as the epimorphism onto
$\mathbf{P}_{k}$. Furthermore, $\mathbf{u}_{k}$ is defined as the composition%
\begin{align}
&  \text{Hom}(\mathbf{K}_{k}\bigcirc\mathbf{K}_{k-1}\bigcirc\cdots
\bigcirc\mathbf{K}_{1},\mathbf{P})\overset{\text{inclusion}}{\hookrightarrow
}\text{Hom}(\mathbf{K}_{k}\otimes\mathbf{K}_{k-1}\bigcirc\cdots\bigcirc
\mathbf{K}_{1},\mathbf{P})\\
&  \cong\text{Hom}(\mathbf{K}_{k},\text{Hom}(\mathbf{K}_{k-1}\bigcirc
\cdots\bigcirc\mathbf{K}_{1},\mathbf{P}))\overset{\underrightarrow
{\text{Hom}(id_{\mathbf{K}_{k}},\mathbf{c}_{k-1})}}{}\text{Hom}(\mathbf{K}%
_{k},\mathbf{P}_{k-1})\nonumber\\
&  \overset{\underrightarrow{\text{Hom}(id_{\mathbf{K}_{k}},\mathbf{e}_{k})}%
}{}\text{Hom}(\mathbf{K}_{k},\mathbf{Q}_{k})\nonumber
\end{align}
([B, Theorem 7.14]).

(8) For a smooth map germ $f:N\rightarrow P$ such that $j^{\infty}f$ is
transverse to $\Sigma^{I_{k}}(N,P)$, let $S^{I_{k}}(j^{\infty}f)$ denote
$(j^{\infty}f)^{-1}(\Sigma^{I_{k}}(N,P))$. If $f|S^{I_{k}}(j^{\infty
}f):S^{I_{k}}(j^{\infty}f)\rightarrow P$ is of kernel rank $i_{k+1}$ at $x$,
then $j_{x}^{\infty}f\in\Sigma^{I_{k+1}}(N,P)$.

(9) The submanifold $\Sigma^{I_{k}}(N,P)$ is actually defined so that it
coincides with the inverse image of the submanifold $\widetilde{\Sigma}%
^{I_{k}}(N,P)$ in $J^{k}(N,P)$ by $\pi_{k}^{\infty}$. The codimension of
$\Sigma^{I_{k}}(N,P)$ in $J^{\infty}(N,P)$ is described in [B, Theorem 6.1].

\begin{remark}
(1) It is known that $\Omega^{I}(N,P)$ is an open subset of $J^{\infty}(N,P)$:
Let $I=(i_{1},i_{2},\cdots,i_{r})$. We prove that the closure of $\Sigma
^{I}(N,P)$ is contained in the subset which consists of all submanifolds
$\Sigma^{J}(N,P)$ of the symbol $J$ of\ length $r$ with $J\geq I$. Let $z\in
J^{\infty}(N,P)$ lies in the closure of $\Sigma^{I}(N,P)$. By definition, we
first have $\dim((\mathrm{Ker}(\mathbf{d}_{1,z}))\geq i_{1}$. If the symbol of
$z$ is $J$ with $J\neq I$, then we can inductively prove that $z$ has a number
$k$ such that $\mathrm{\dim}(\mathrm{Ker}(\mathbf{d}_{j,z}|\mathbf{K}%
_{j-1,z}))=i_{j}$ for $1\leq j\leq k<r$ and $\dim(\mathrm{Ker}(\mathbf{d}%
_{k+1,z}|\mathbf{K}_{k,z}))>i_{k+1}$. This implies the assertion.

(2) If a symbol $J$ is an infinite series $(j_{1},j_{2},\cdots,j_{k},\cdots)$
and $\mathrm{codim}\Sigma^{J}(N,P)\leq n$, then $j_{1},j_{2},\cdots
,j_{k},\cdots$ are equal to $0$ except for a finite number of $j_{k}$'s.
\end{remark}

\section{Polynomials}

Let $V$ and $\ W$ be vector spaces with inner product of dimensions $v$ and
$w$ respectively. Let $e_{1},e_{2},\cdots,e_{v}$ and $d_{1},d_{2},\cdots
,d_{w}$ be orthogonal basis of $V$ and $W$ respectively. We introduce the
inner product in Hom$(\otimes^{\ell}V,W)$ as follows. Let $h_{i}%
\in\mathrm{Hom}(\otimes^{\ell}V,W)$ ($i=1,2$) and let%
\[
h_{1}(e_{i_{1}}\otimes\ldots\otimes e_{i_{\ell}})=\sum_{j=1}^{w}a_{i_{1}%
i_{2}\cdots i_{\ell}}^{j}d_{j}\text{ \ \ and \ \ }h_{2}(e_{i_{1}}\otimes
\ldots\otimes e_{i_{\ell}})=\sum_{j=1}^{w}b_{i_{1}i_{2}\cdots i_{\ell}}%
^{j}d_{j}.
\]
Then we define the inner product by%
\[
\langle h_{1},h_{2}\rangle=\sum_{j=1}^{w}\sum_{i_{1}i_{2}\cdots i_{\ell}%
}a_{i_{1}i_{2}\cdots i_{\ell}}^{j}b_{i_{1}i_{2}\cdots i_{\ell}}^{j}.
\]
Let $S$ and $T$ be isomorphisms of $V$ and $W$ which preserve the inner
products respectively.\ We define the action of $(T,S)$ on Hom$(\otimes^{\ell
}V,W)$ by $(T,S)h=T\circ h\circ(\otimes^{\ell}S^{-1})$. We show by induction
on $\ell$\ that this inner product is invariant with respect to this action.
We represent $S^{-1}$ by the matrix $(s_{ij})$ under the basis $e_{1}%
,e_{2},\cdots,e_{w}$.

The assertion for $\ell=1$ is well known. Assume that the assertion holds for
$\ell-1$. Under the canonical isomorphism Hom$(\otimes^{\ell}V,W)\cong
\mathrm{Hom}(V,\mathrm{Hom}(\otimes^{\ell-1}V,W))$ we let $h\in\mathrm{Hom}%
(\otimes^{\ell}V,W)$ correspond to $\overline{h}$, which satisfies
$\overline{h}(e_{i_{1}})(e_{i_{2}}\otimes\ldots\otimes e_{i_{\ell}%
})=h(e_{i_{1}}\otimes e_{i_{2}}\otimes\ldots\otimes e_{i_{\ell}})$. Then we
have that $\langle h_{1},h_{2}\rangle=\Sigma_{j=1}^{v}\langle\overline{h_{1}%
}(e_{j}),\overline{h_{2}}(e_{j})\rangle$. Hence, we have that
\begin{align*}
\langle(T,S)h_{1},(T,S)h_{2}\rangle &  =\sum_{i=1}^{v}\langle((T,S)\overline
{h_{1}})(S^{-1}(e_{i})),((T,S)\overline{h_{2}})(S^{-1}(e_{i}))\rangle\\
&  =\sum_{i=1}^{v}\langle\overline{h_{1}}(S^{-1}(e_{i})),\overline{h_{2}%
}(S^{-1}(e_{i}))\rangle\\
&  =\sum_{i=1}^{v}\langle\overline{h_{1}}(\Sigma_{j=1}^{v}s_{ij}%
e_{j}),\overline{h_{2}}(\Sigma_{k=1}^{v}s_{ik}e_{k})\rangle\\
&  =\sum_{i=1}^{v}(\Sigma_{j=1}^{v}s_{ij}(\Sigma_{k=1}^{v}s_{ik}%
\langle\overline{h_{1}}(e_{j}),\overline{h_{2}}(e_{k})))\rangle\\
&  =\sum_{j=1}^{v}\sum_{k=1}^{v}(\Sigma_{i=1}^{v}s_{ij}s_{ik})\langle
\overline{h_{1}}(e_{j}),\overline{h_{2}}(e_{k})\rangle\\
&  =\sum_{j=1}^{v}\sum_{k=1}^{v}\delta_{jk}\langle\overline{h_{1}}%
(e_{j}),\overline{h_{2}}(e_{k})\rangle\\
&  =\sum_{j=1}^{v}\langle\overline{h_{1}}(e_{j}),\overline{h_{2}}%
(e_{j})\rangle\\
&  =\langle h_{1},h_{2}\rangle.
\end{align*}

We recall that Hom$(\Sigma_{j=1}^{\ell}\bigcirc^{j}V,W)$ is identified with
the set of polynomials of degree $\leq\ell$ having the constant $0$ (see
[Mats, Ch. III, Section 2]). Let $\mathbf{V}$ and $\mathbf{W}$ be smooth
vector bundles with metric over a manifold $S$ with fibers $V$ and $W$
respectively. Then Hom$(\bigcirc^{\ell}\mathbf{V},\mathbf{W})$ is also a
vector bundle with metric. For a point $c\in S$, take an open neighborhood $U$
around $c$ such that $\mathbf{V}|_{U}$ and $\mathbf{W}|_{U}$ are the trivial
bundles, say $U\times V$ and $U\times W$ respectively. Then an element of
Hom$(\bigcirc^{\ell}\mathbf{V},\mathbf{W})|_{U}$ is identified with a
polynomial $\Sigma_{j=1}^{w}(\Sigma_{|\omega|=\ell}A_{\omega}^{j}%
(c)x_{1}^{\omega_{1}}x_{2}^{\omega_{2}}\cdots x_{v}^{\omega_{v}})d_{j}$, $c\in
U$, where $\omega=(\omega_{1},\omega_{2},\cdots,\omega_{v})$, $\omega_{i}%
\geq0$ ($i=1,\cdots,v$), and $|\omega|=\omega_{1}+\cdots+\omega_{v}$ and
$A_{\omega}^{j}(c)$ is a real number. If $A_{\omega}^{j}(c)$ are smooth
functions of $c,$ then $\{A_{\omega}^{j}(c)\}$ defines a smooth section of
Hom$(\bigcirc^{\ell}\mathbf{V},\mathbf{W})|_{U}$ over $U$.

We now provide $N$ and $P$ with Riemannian metrics respectively. Then they
induce the metrics on $\mathbf{D}$ and $\mathbf{P}$, and hence induces the
metric on Hom$(\mathbf{K}_{k}\bigcirc\mathbf{K}_{k-1}\bigcirc\cdots
\bigcirc\mathbf{K}_{1},\mathbf{P})$. Furthermore, we can prove inductively
that $\mathbf{P}_{k},$ and also $\mathbf{Q}_{k+1}$ as the orthogonal
complement of Im$(d_{k+1}|\mathbf{K}_{k})$\ inherit the induced metrics by (5)
and (6) in Section 2 respectively. Consequently we have the induced metric on
Hom$(\mathbf{K}_{k+1},$ $\mathbf{Q}_{k+1})$.

Let us recall $\mathbf{d}_{k}|\mathbf{K}_{k-1}:\mathbf{K}_{k-1}\rightarrow
\mathbf{P}_{k-1}$ and $\mathbf{e}_{k-1}:\mathbf{P}_{k-1}\rightarrow
\mathbf{Q}_{k}$ over $\Sigma^{I_{k}}(N,P),$ which induces the commutative
diagram%
\[%
\begin{array}
[c]{lllllll}%
\mathbf{K}_{k} & \rightarrow & \mathbf{K}_{k-1} & \rightarrow & \text{
\ \ \ \ \ \ \ }\mathbf{P}_{k-1} & \rightarrow & \mathbf{Q}_{k}\\
&  & \downarrow &  & \text{ \ \ \ \ \ \ \ }\curvearrowright & \swarrow
\mathbf{j}_{\mathbf{Q}_{k}} & \\
\mathbf{0} & \rightarrow & \mathbf{K}_{k-1}/\mathbf{K}_{k} & \rightarrow &
\text{Hom}(\mathbf{K}_{k-1},\mathbf{Q}_{k-1}). &  &
\end{array}
\]
Since $\mathbf{Q}_{k}$ is the cokernel of $\mathbf{d}_{k}|\mathbf{K}_{k-1}$,
we obtain the canonical isomorphism%
\[
\mathbf{j}_{\mathbf{Q}_{k}}:\mathbf{Q}_{k}\rightarrow\mathrm{\operatorname{Im}%
}(\mathbf{d}_{k}|\mathbf{K}_{k-1})^{\bot}\text{ \ \ \ \ over \ }\Sigma^{I_{k}%
}(N,P),
\]
where the symbol $\bot$ refers to the orthogonal complement. We also use the
notation\ $\mathbf{j}_{\mathbf{Q}_{k}}:\mathbf{Q}_{k}\rightarrow
\mathrm{Hom}(\mathbf{K}_{k-1},\mathbf{Q}_{k-1})$.

Let $k\geq2$. We now construct the homomorphism, for $1\leq i\leq k$,%
\begin{equation}
\mathbf{q}(k)^{i+1,i+1}:T(\Sigma^{I_{i-1}}(N,P))|_{\Sigma^{I_{k}}%
(N,P)}\bigcirc\mathbf{K}_{i}\bigcirc\mathbf{K}_{i-1}\bigcirc\cdots
\bigcirc\mathbf{K}_{1}\rightarrow\mathbf{Q}_{1}%
\end{equation}
over $\Sigma^{I_{k}}(N,P)$ inductively by using $\mathbf{d}_{i+1}%
|_{\Sigma^{I_{k}}(N,P)}:T(\Sigma^{I_{i-1}}(N,P))|_{\Sigma^{I_{k}}%
(N,P)}\rightarrow\mathbf{P}_{i}|_{\Sigma^{I_{k}}(N,P)}$\ as follows. By the
inclusion $\mathbf{P}_{i}|_{\Sigma^{I_{k}}(N,P)}\subset\mathrm{Hom}%
(\mathbf{K}_{i},\mathbf{Q}_{i})|_{\Sigma^{I_{k}}(N,P)}$ we have the
homomorphism%
\[
\mathbf{q}(k)_{\otimes}^{i+1,2}:(T(\Sigma^{I_{i-1}}(N,P))|_{\Sigma^{I_{k}%
}(N,P)})\otimes\mathbf{K}_{i}\rightarrow\mathbf{Q}_{i}\text{ \ \ \ \ over
}\Sigma^{I_{k}}(N,P).
\]
Suppose that we have constructed the homomorphism, for $j\leq i$,%
\[
\mathbf{q}(k)_{\otimes}^{i+1,i-j+2}:T(\Sigma^{I_{i-1}}(N,P))|_{\Sigma^{I_{k}%
}(N,P)}\otimes\mathbf{K}_{i}\otimes\mathbf{K}_{i-1}\otimes\cdots
\otimes\mathbf{K}_{j}\rightarrow\mathbf{Q}_{j}%
\]
over $\Sigma^{I_{k}}(N,P)$. By using $\mathbf{j}_{\mathbf{Q}_{j}}%
:\mathbf{Q}_{j}\rightarrow$Hom$(\mathbf{K}_{j-1},\mathbf{Q}_{j-1})$\ over
$\Sigma^{I_{k}}(N,P)$, we obtain the homomorphism%
\begin{equation}
\mathbf{q}(k)_{\otimes}^{i+1,i-j+3}:T(\Sigma^{I_{i-1}}(N,P))|_{\Sigma^{I_{k}%
}(N,P)}\otimes\mathbf{K}_{i}\otimes\mathbf{K}_{i-1}\otimes\cdots
\otimes\mathbf{K}_{j-1}\rightarrow\mathbf{Q}_{j-1}%
\end{equation}
over $\Sigma^{I_{k}}(N,P)$. By setting $j=2$, we obtain $\mathbf{q}%
(k)_{\otimes}^{i+1,i+1}$. It remains to prove that $\mathbf{q}(k)_{\otimes
}^{i+1,i+1}$ is symmetric. This fact has been esssentially stated in [B,
Section 7, p.413] without proof. Following the proof of [B, Theorem 4.1] we
briefly prove it.

Let $z\in\Sigma^{I_{k}}(N,P)$. By the Riemannian metric of $P$, we consider
the convex neighborhood of $P$ around $\pi_{P}^{\infty}(z)=y$. Let us
canonically identify $\mathbf{Q}_{1,z}$ with a subspace of $T_{y}P$ by the
isomorphism $\mathbf{P}_{z}\rightarrow T_{y}P$. By taking a basis of
$\mathbf{Q}_{1,z}$ and projecting it by the exponential map, we have the local
coordinates $y_{1},y_{2},\cdots,y_{p-n+i_{1}}$ on the convex neighborhood of
$y$. Then we identify $\mathbf{Q}_{1,z}$ with Hom$(\frak{m}_{y}^{\mathbf{Q}%
}/(\frak{m}_{y}^{\mathbf{Q}})^{2},\mathbf{R})$, where $\frak{m}_{y}%
^{\mathbf{Q}}$ is the ideal generated by $y_{1},y_{2},\cdots,y_{p-n+i_{1}}$.
Let $D$ and $D_{j}$ be sections of $T(\Sigma^{I_{i-1}}(N,P))|_{\Sigma^{I_{k}%
}(N,P)}$ and $\mathbf{K}_{j}$ defined around $z$ and let $\alpha\in
\frak{m}_{y}^{\mathbf{Q}}$. Then (3.1) is regarded as the homomorphism induced
from%
\[
T(\Sigma^{I_{i-1}}(N,P))_{z}\otimes\mathbf{K}_{i,z}\otimes\mathbf{K}%
_{i-1,z}\otimes\cdots\otimes\mathbf{K}_{1,z}\otimes\frak{m}_{y}^{\mathbf{Q}%
}/(\frak{m}_{y}^{\mathbf{Q}})^{2}\longrightarrow\mathbf{R}%
\]
which maps $D\otimes D_{i}\otimes\cdots\otimes D_{1}\otimes\alpha$ to
$(DD_{i}\cdots D_{1}\alpha)(z)$ (see (a) and (b) in the proof of [B, Theorem
4.1]). We have to show the following for the symmetry (consult Remark 3.1
below to avoid the infinity of the dimensions of the tangent spaces). In the
expression with $[D_{j},D_{j-1}]=D_{j}D_{j-1}-D_{j-1}D_{j}$%
\[
DD_{i}\cdots D_{j}D_{j-1}\cdots D_{1}\alpha-DD_{i}\cdots D_{j-1}D_{j}\cdots
D_{1}\alpha=DD_{i}\cdots\lbrack D_{j},D_{j-1}]\cdots D_{1}\alpha
\]
for some $j$ with $1<j\leq i+1$ ($D_{i+1}=D$), we have that $[D_{j},D_{j-1}]$
is the section of $\mathbf{K}_{j-1}$ for $j\leq i$ and of $T(\Sigma^{I_{i-1}%
}(N,P))_{z}$\ for $j=i+1$\ by [B, Lemma 3.2]. Since $\mathbf{K}_{j}%
|_{\Sigma^{I_{k}}(N,P)}\subset\mathbf{K}_{j-1}|_{\Sigma^{I_{k}}(N,P)}$, the
length of $DD_{i}\cdots\lbrack D_{j},D_{j-1}]\cdots D_{1}$ is $i$, $D$ and
$\ [D,D_{i}]$ lie in $T(\Sigma^{I_{i-1}}(N,P))_{z}$ and since $T(\Sigma
^{I_{i-1}}(N,P))_{z}\subset T(\Sigma^{I_{i-2}}(N,P))_{z}$, we have that
$(DD_{i}\cdots\lbrack D_{j},D_{j-1}]\cdots D_{1}\alpha)(z)=0$ by
Ker$(\mathbf{d}_{i,z})=T(\Sigma^{I_{i-1}}(N,P))_{z}$ in (2.1). This is what we want.

In particular, if $i=1$ and we restrict $T(\Sigma^{I_{i-1}}(N,P))|_{\Sigma
^{I_{k}}(N,P)}$ to $\mathbf{K}_{1},$ then we have the homomorphism
$\mathbf{q}(k)^{2,2}|(\mathbf{K}_{1}\bigcirc\mathbf{K}_{1}):\mathbf{K}%
_{1}\bigcirc\mathbf{K}_{1}\rightarrow\mathbf{Q}_{1}$\ over $\Sigma^{I_{k}%
}(N,P)$, which induces the nonsingular quadratic form $(\mathbf{K}%
_{1}/\mathbf{K}_{2})\bigcirc(\mathbf{K}_{1}/\mathbf{K}_{2})\rightarrow
\mathbf{Q}_{1}$ on each fiber.

\begin{remark}
We can entirely do the arguments in Sections 2 and 3 on $J^{\ell}(N,P)$ for a
large $\ell$. We provide $N$ and $P$ with Riemannian metrics. For any points
$x\in N$ and $y\in P$, we have the local coordinates $(x_{1},...,x_{n})$ and
$(y_{1},...,y_{n})$ on convex neighbourhoods of $x$ and $y$ associated to
orthonormal basis of $T_{x}N$ and $T_{y}P$ respectively. Let us define the
canonical embedding $\mu_{\infty}^{\ell}:J^{\ell}(TN,TP)\rightarrow J^{\infty
}(TN,TP)$ such that $\pi_{\ell}^{\infty}\circ\mu_{\infty}^{\ell}=id_{J^{\ell
}(TN,TP)}$ and that the $i$-th components for $i>\ell$ of elements of the
image $\mu_{\infty}^{\ell}$ are the null homomorphisms of $\mathrm{Hom}%
(S^{i}(\pi_{N}^{\ast}(TN)),\pi_{P}^{\ast}(TP))$. We regard $\mu_{\infty}%
^{\ell}$\ as the map to $J^{\infty}(N,P)$\ under the identification (1.2). Any
element $z\in\mu_{\infty}^{\ell}(J^{\ell}(TN,TP))$ is represented by a
$C^{\infty}$ map germ $f:(N,x)\rightarrow(P,y)$ such that any $i$-th
derivative of $f$ with $i>\ell$ vanishes under these coordinates. It is clear
that We can prove that $\mathbf{D}|_{\mu_{\infty}^{\ell}(J^{\ell}(TN,TP))}$ is
tangent to $\mu_{\infty}^{\ell}(J^{\ell}(TN,TP))$. Indeed, for $\sigma
=(\sigma_{1},...,\sigma_{n})$ with non-negative integers $\sigma_{i}$, let us
recall the functions $X_{i}$ and $Z_{j,\sigma}$ with $1\leq i\leq n$ and
$1\leq j\leq p$ defined locally on a neighbourhood of $J^{\infty}(N,P)$ by,
for $z=j_{x}^{\infty}f,$
\[
X_{i}(z)=x_{i}\text{ \ and \ }Z_{j,\sigma}(z)=\frac{\partial^{|\sigma|}%
(y_{j}\circ f)}{\partial x_{1}^{\sigma_{1}}\cdots\partial x_{n}^{\sigma^{n}}%
}(x),
\]
which constitute the local coordinates on $J^{\infty}(N,P)$ as described in
[B, Section 1]. Let $\Phi$ be a smooth function defined locally on
$\mu_{\infty}^{\ell}(J^{\ell}(TN,TP))$ and let $D_{i}\in\mathbf{D}$ be the
total tangent vector corresponding to $\partial/\partial x_{i}$ by the
canonical identification of $\mathbf{D}$ and $(\pi_{N}^{\infty})^{\ast}(TN)$.
Let $\sigma^{\prime}=(\sigma_{1},...,\sigma_{i-1,}\sigma_{i}+1,\sigma
_{i+1},...,\sigma_{n})$. Then we have
\[
D_{i}(\Phi)(z)=\frac{\partial(\Phi\circ j^{\infty}f)}{\partial x_{i}}%
(x)=\frac{\partial\Phi}{\partial X_{i}}(z)+\underset{j,\sigma}{\sum}%
\frac{\partial\Phi}{\partial Z_{j,\sigma}}(z)Z_{j,\sigma^{\prime}}(z)
\]
by [B, (1.8)]. If $z\in\mu_{\infty}^{\ell}(J^{\ell}(TN,TP))$, then
$Z_{j,\sigma}(z)$ vanishes for $|\sigma|>\ell$. Hence, $D_{i}(\Phi)$ is a
smooth function defined locally on $\mu_{\infty}^{\ell}(J^{\ell}(TN,TP))$.
This implies that $D_{i}$ is tangent to $\mu_{\infty}^{\ell}(J^{\ell}%
(TN,TP))$. Since $\mathbf{D}_{z}$ consists of all linear combinations of
$D_{1},...,D_{n}$, we have that $\mathbf{D}_{z}\subset T_{z}(\mu_{\infty
}^{\ell}(J^{\ell}(TN,TP)))$. Therefore, we can do the required arguments on
$\mu_{\infty}^{\ell}(J^{\ell}(TN,TP))$, namely also on $J^{\ell}(TN,TP)$.
\end{remark}

\section{Primary obstruction}

Let $\Gamma_{\Omega^{I}}^{tr}(N,P)$ denote the subspace of $\Gamma_{\Omega
^{I}}(N,P)$ consisting of all continuous sections of $\pi_{N}^{\infty}%
|\Omega^{I}(N,P):\Omega^{I}(N,P)\rightarrow N$ which are transverse to each
$\Sigma^{J}(N,P)$. For $\frak{s}\in\Gamma_{\Omega^{I}}^{tr}(N,P)$, we set
$S^{I_{j}}(\frak{s})=\frak{s}^{-1}(\Sigma^{I_{j}}(N,P))$, $S^{I_{j}%
,0}(\frak{s})=\frak{s}^{-1}(\Sigma^{I_{j},0}(N,P))$, $(\frak{s}|S^{I}%
(\frak{s}))^{\ast}(\mathbf{K}_{j})=K_{j}(S^{I}(\frak{s}))$ and $(\frak{s}%
|S^{I}(\frak{s}))^{-1}\mathbf{Q}_{1}=Q(S^{I}(\frak{s}))$. We often write
$S^{I}(\frak{s})$ as $S^{I}$ if there is no confusion.

Let $L=(\ell_{1},\ell_{2},\cdots,\ell_{r})$ and\ $I=(i_{1},i_{2},\cdots
,i_{k},0,\cdots,0,\cdots)$ such that $I_{r}\leq L$, codim $\Sigma^{L}(n,p)\leq
n$ and codim $\Sigma^{I}(n,p)\leq n$, where $k$ may be larger than $r$. Let
$C(I^{+})$ (resp. $C(I)$) refer to the union $C\cup(\cup_{J>I}S^{J}(s))$
(resp. $C\cup(\cup_{J\geq I}S^{J}(s))$), where $J$ are symbols of infinite
length and $C$ is a closed subset of $N$.

We show in this section that it is enough for the proof of Theorem 0.1 to
prove Theorem 4.1.

\begin{theorem}
Let $n\geq p\geq2$. Let $N$ and $P$ be connected manifolds of dimensions $n$
and $p$ respectively with $\partial N=\emptyset$. Let $L$ and $I=(i_{1}%
,i_{2},\cdots,i_{k},0)$ be as above. Assume that $\Omega^{L}(N,P)$ contains
$\Sigma^{n-p+1,0}(N,P)$ at least. Let $s$ be a section of $\Gamma_{\Omega^{L}%
}^{tr}(N,P)$ which has an $\Omega^{L}$-regular map $g(I^{+})$ defined on a
neighborhood of $C(I^{+})$ into $P$, where $j^{\infty}g(I^{+})=s$. Then there
exists a homotopy $s_{\lambda}\in\Gamma_{\Omega^{L}}^{tr}(N,P)$ relative to a
neighborhood of $C(I^{+})$ with the following properties.

$(1)$ $s_{0}=s$ and $s_{1}\in\Gamma_{\Omega^{L}}^{tr}(N,P)$.

$(2)$ There exists an $\Omega^{L}$-regular map $g_{I}$ defined on a
neighborhood of $C(I)$, where $j^{\infty}g_{I}=s_{1}$ holds.

$(3)$ $s^{-1}(\Sigma^{I}(N,P))=(j^{\infty}g_{I})^{-1}(\Sigma^{I}(N,P))$.
\end{theorem}

The case $I=(n-p+1,0)$ of Theorem 4.1 follows from Theorem 1 of [An1], where a
partially sketchy proof was given and the detailed proof was given in [An4,
Theorem 4.1] and [An6, Theorem 0.5]. Let us explain how it follows. In fact,
we have $\Omega^{n-p+1,0}(N,P)=\Sigma^{n-p}(N,P)\cup\Sigma^{n-p+1,0}(N,P)$,
and if we set $N_{0}=S^{n-p}(s)\cup S^{n-p+1,0}(s)$, then $C((n-p+1,0)^{+}%
)=C\cup(N\setminus N_{0})$. Let $U$ and $U^{\prime}$ be closed neighborhoods
of $C\cup(N\setminus N_{0})$ with $U\subset$Int$U^{\prime}$, where
$g((n-p+1,0)^{+})$ is defined. Since $s\in\Gamma_{\Omega^{L}}^{tr}(N,P)$,
$s|N_{0}$ is a section of $\Gamma_{\Omega^{n-p+1,0}}^{tr}(N_{0},P)$ and
$g((n-p+1,0)^{+})|(U^{\prime}\cap N_{0})$ is an $\Omega^{n-p+1,0}$-regular
map. Hence, we obtain a homotopy $u_{\lambda}\in\Gamma_{\Omega^{n-p+1,0}}%
^{tr}(N_{0},P)$ relative to a neighborhood of $U\cap N_{0}$ and an
$\Omega^{n-p+1,0}$-regular map $f_{0}:N_{0}\rightarrow P$ such that
$s_{0}|N_{0}=u_{0}$ and $u_{1}=j^{\infty}f_{0}$. Then we obtain a required
homotopy $s_{\lambda}$ by defining $s_{\lambda}|N_{0}=u_{\lambda}$ and
$s_{\lambda}|U=j^{\infty}g((n-p+1,0)^{+})$.

We will prove Theorem 4.1 for $I>(n-p+1,0)$ in Section 6.

We now prove Theorem 0.1 for $\Omega^{L}$ for this symbol $L=(\ell_{1}%
,\ell_{2},\cdots,\ell_{r})$ in place of $\Omega^{I}$\ by using Theorem 4.1. In
Sections 4, 5 and 6 we use the notation $\Omega$ for $\Omega^{L}$.

\begin{proof}
[Proof of Theorem 0.1]Suppose that the section $s$ given in Theorem 0.1 lies
in $\Gamma_{\Omega}^{tr}(N,P)$. Let $I=(i_{1},i_{2},\cdots,i_{k}%
,0,\cdots,0,\cdots)$ be the largest symbol such that $I_{\ell}\leq L$ and
codim $\Sigma^{I}(n,p)\leq n$. We can choose such a symbol $I$ by using
Section 2 (4) and Remark 2.1 (2). Then we first set $C(I^{+})=C$ and
$g(I^{+})=g$. By Theorem 4.1 there exists an $\Omega$-regular map $g_{I}$
defined on a neighborhood of $C(I)$, where $j^{\infty}g_{I}=s$ holds. If we
note Remark 2.1 (2), then we can prove Theorem 0.1 by the downward induction
on the symbols $I$ in the lexicographic order.
\end{proof}

We begin by preparing several notions and results, which are necessary for the
proof of Theorem 4.1. For the map $g(I^{+})$ and the closed subset $C(I^{+})$,
we take an open neighborhood $U(C(I^{+}))^{\prime}$ of $C(I^{+})$, where
$j^{\infty}g(I^{+})=s$. Without loss of generality we may assume that
$N\setminus U(C(I^{+}))^{\prime}$ is nonempty. Take a smooth function
$h_{C(I^{+})}:N\rightarrow\lbrack0,1]$ such that
\begin{equation}
\left\{
\begin{array}
[c]{ll}%
h_{C(I^{+})}(x)=1 & \text{for }x\in C(I^{+}),\\
h_{C(I^{+})}(x)=0 & \text{for }x\in N\setminus U(C(I^{+}))^{\prime},\\
0<h_{C(I^{+})}(x)<1 & \text{for }x\in U(C(I^{+}))^{\prime}\setminus C(I^{+}).
\end{array}
\right.
\end{equation}
By the Sard Theorem ([H2]) there is a regular value $r$ of $h_{C(I^{+})}$ with
$0<r<1$. Then $h_{C(I^{+})}^{-1}(r)$ is a submanifold and we set
$U(C(I^{+}))=h_{C(I^{+})}^{-1}([r,1])$. We decompose $N\setminus
\mathrm{Int}U(C(I^{+}))$ into the connected components, say $L_{1}%
,\cdots,L_{j},\cdots$. It suffices to prove Theorem 4.1 for each $L_{j}\cup
$Int$U(C(I^{+}))$. Since $\partial N=\emptyset$, we have that $N\setminus
U(C(I^{+}))$ has empty boundary. If $L_{j}$ is not compact, then Theorem 4.1
holds for $L_{j}\cup$Int$U(C(I^{+}))$ by Gromov's theorem ([G1, Theorem
4.1.1]). Therefore, it suffices to consider the special case where

(C1) $N\setminus\mathrm{Int}U(C(I^{+}))$ is compact, connected and nonempty,

(C2) $\partial U(C(I^{+}))$ is a submanifold of dimension $n-1$,

(C3) for the smooth function $h_{C(I^{+})}:N\rightarrow\lbrack0,1]$ satisfying
(4.1) there is a sufficiently small positive real number $\varepsilon$ with
$r-2\varepsilon>0$ such that $r-t\varepsilon$ ($0\leq t\leq2$) are all regular
values of $h_{C(I^{+})}$. We have that $h_{C(I^{+})}^{-1}([r-2\varepsilon,1])$
is contained in $U(C(I^{+}))^{\prime}$.

We set $U(C(I^{+}))_{t}=h_{C(I^{+})}^{-1}([r-(2-t)\varepsilon,1])$. In
particular, we have $U(C(I^{+}))_{2}=U(C(I^{+}))$. Furthermore, we may assume that

(C4) $s\in\Gamma_{\Omega}^{tr}(N,P)$ and $S^{I}(s)$ is transverse to $\partial
U(C(I^{+}))_{0}$ and $\partial U(C(I^{+}))_{2}$.

In what follows we choose and fix a Riemannian metric of $N$, which satisfies

\textbf{Orthogonality Condition}; \textit{for the symbol }$I$\textit{,
}$K_{j-1}(S^{I}(s))/K_{j}(S^{I}(s))$\textit{ is orthogonal to }$S^{I_{j}}%
(s)$\textit{ in }$S^{I_{j-1}}(s)$\textit{ for }$k\leq j\leq1$\textit{ on
}$S^{I}(s)$\textit{ }$(S^{I_{0}}(s)=N)$.

Let $\nu(\Sigma^{I})$ be the normal bundle $(T(J^{\infty}(N,P))|_{\Sigma^{I}%
})/T(\Sigma^{I}(N,P))$ and let $c(I)=\dim\nu(\Sigma^{I})$. Let us fix a direct
sum decomposition%
\begin{equation}
\nu(\Sigma^{I}){{{{=\oplus}}}}_{j=1}^{k}\nu(I_{j}\subset I_{j-1})|_{\Sigma
^{I}(N,P)},
\end{equation}
and the direct sum decomposition $\mathbf{K}_{1}={{{{\oplus}}}}_{j=1}%
^{k-1}(\mathbf{K}_{j}\mathbf{/K}_{j+1})\oplus\mathbf{K}_{k}$ over $\Sigma
^{I}(N,P)$. Let $\mathbf{j}_{\mathbf{K}}:\mathbf{K}_{1}\rightarrow\nu
(\Sigma^{I})$ over $\Sigma^{I}(N,P)$ be the composition of the inclusion
$\mathbf{K}_{1}\rightarrow T(J^{\infty}(N,P))$ and the projection
$T(J^{\infty}(N,P))|_{\Sigma^{I}(N,P)}\rightarrow\nu(\Sigma^{I})$. We have the
monomorphism%
\[
\mathbf{j}_{\mathbf{K}}\circ(s|S^{I})^{\mathbf{K}_{1}}:K_{1}(S^{I}%
(s))\rightarrow\mathbf{K}_{1}|_{\Sigma^{I}(N,P)}\rightarrow\nu(\Sigma^{I}).
\]

For $s\in\Gamma_{\Omega}^{tr}(N,P)$, let $\frak{n}(s,I)$ or simply
$\frak{n}(I)$ be the orthogonal normal bundle of $S^{I}(s)$ in $N$. Let
$\frak{n}(s,I_{j}\subset I_{j-1})$ be the orthogonal normal bundle of
$S^{I_{j}}(s)$ in $S^{I_{j-1}}(s)$ over $S^{I}(s)$. Then we have the canonical
direct sum decomposition such as%
\begin{equation}
\frak{n}(s,{{{{I)=\oplus}}}}_{j=1}^{k}\frak{n}(s,I_{j}\subset I_{j-1}),
\end{equation}
Furthermore, we obtain the bundle map
\[
ds|\frak{n}(s,I):\frak{n}(s,I)\rightarrow\nu(\Sigma^{I})
\]
covering $s|S^{I}:S^{I}(s)\rightarrow\Sigma^{I}(N,P)$. Let $\mathbf{i}%
_{\frak{n}(s,I)}:\frak{n}(s,I)\subset TN|_{S^{I}}$ denote the inclusion. We
define $\Psi(s,I):K_{1}(S^{I}(s))\rightarrow\frak{n}(s,I)\subset TN|_{S^{I}}$
to be the composition%
\begin{align*}
&  \mathbf{i}_{\frak{n}(s,I)}\circ((s|S^{I})^{\ast}(ds|\frak{n}(s,I)))^{-1}%
\circ((s|S^{I})^{\ast}(\mathbf{j}_{\mathbf{K}}\circ(s|S^{I})^{\mathbf{K}_{1}%
}))\\
\text{ \ \ \ \ \ \ }  &  :K_{1}(S^{I}(s))\rightarrow(s|S^{I})^{\ast}\nu
(\Sigma^{I})\rightarrow\frak{n}(s,I)\hookrightarrow TN|_{S^{I}}.
\end{align*}
We note that this homomorphism does not use the decomposition in (4.2) and we
can take the direct sum decompositions in (4.2) to be compatible with those in
(4.3). Let $i_{K_{1}(S^{I}(s))}:K_{1}(S^{I}(s))\rightarrow TN|_{S^{I}}$ be the inclusion.

\begin{remark}
If $f$ is an $\Omega$-regular map, then it follows from the definition of
$\mathbf{D}$ that $\mathbf{i}_{K_{1}(S^{I}(j^{\infty}f))}=\Psi(j^{\infty}f,I)$.
\end{remark}

In what follows let $M=S^{I}(s)\setminus$Int$(U(C(I^{+})))$. Let
Mono$(K_{1}(S^{I}(s))|_{M},TN|_{M})$ denote the subset of Hom$(K_{1}%
(S^{I}(s))|_{M},TN|_{M})$ which consists of all monomorphisms $K_{1}%
(S^{I}(s))_{c}\rightarrow T_{c}N$, $c\in M$. We denote the bundle of the local
coefficients $\mathcal{B}(\pi_{j}(\mathrm{Mono}(K_{1}(S^{I}(s))_{c}%
,T_{c}N))),$ $c\in M,$ by $\mathcal{B}(\pi_{j})$, which is a covering space
over $M$ with fiber $\pi_{j}(\mathrm{Mono}(K_{1}(S^{I}(s))_{c},T_{c}N))$
defined in [Ste, 30.1]. From the obstruction theory due to [Ste, 36.3], it
follows that the obstructions for $i_{K_{1}(S^{I}(s))}|_{M}$ and
$\Psi(s,I)|_{M}$ to be homotopic are the primary differences $d(i_{K_{1}%
(S^{I}(s))}|_{M},\Psi(s,I)|_{M})$, which are defined in $H^{j}(M,\partial
M;\mathcal{B}(\pi_{j}))$ with the local coefficients . We show that all of
them vanish. In fact, note $I>(n-p+1,0)$. If $i_{1}=n-p+1$, then we have%
\[
\dim M<\dim S^{i_{1}}=n-i_{1}(p-n+i_{1})=n-i_{1}.
\]
If $i_{1}>n-p+1$, then
\[
\dim M\leq\dim S^{i_{1}}=n-i_{1}(p-n+i_{1})<n-i_{1}.
\]
Since \textrm{Mono}$(\mathbf{R}^{i_{1}},\mathbf{R}^{n})$ is identified with
$GL(n)/GL(n-i_{1})$, it follows from [Ste, 25.6 \ 38.2] that $\pi_{j}%
($\textrm{Mono}$(\mathbf{R}^{i_{1}},\mathbf{R}^{n}))\cong\{\mathbf{0}\}$ for
$j<n-i_{1}(\leq p-1)$. Hence, there exists a homotopy $\psi^{M}(s,I)_{\lambda
}:K_{1}(S^{I}(s))|_{M}\rightarrow TN|_{M}$ relative to $M\cap U(C(I^{+}))_{1}$
in \textrm{Mono}$(K_{1}(S^{I}(s))|_{M},TN|_{M})$ such that $\psi^{M}%
(s,I)_{0}=i_{K_{1}(S^{I}(s))}|_{M}$ and $\psi^{M}(s,I)_{1}$$=\Psi(s,I)|_{M}$.
Let $\mathrm{Iso}(TN|_{M},TN|_{M})$\ denote the subspace of $\mathrm{Hom}%
(TN|_{M},TN|_{M})$ which consists of all isomorphisms of $T_{c}N$, $c\in
M$.\ The restriction map
\[
r_{M}:\mathrm{Iso}(TN|_{M},TN|_{M})\rightarrow\mathrm{Mono}(K_{1}%
(S^{I}(s))|_{M},TN|_{M})
\]
defined by $r_{M}(h)=h|(K_{1}(S^{I}(s))_{c})$, for $h\in\mathrm{Iso}%
(T_{c}N,T_{c}N)$, induces a structure of a fiber bundle with fiber
Iso$(\mathbf{R}^{n-i_{1}},\mathbf{R}^{n-i_{1}})\times{\mathrm{Hom}}%
(\mathbf{R}^{n-i_{1}},\mathbf{R}^{i_{1}})$. By applying the covering homotopy
property of the fiber bundle $r_{M}$ to the sections $id_{TN|_{M}}$ and the
homotopy $\psi^{M}(s,I)_{\lambda},$ we obtain a homotopy $\Psi^{M}%
(s,I)_{\lambda}:$ $TN|_{M}\rightarrow TN|_{M}$ such that $\Psi^{M}%
(s,I)_{0}=id_{TN|_{M}}$, $\Psi^{M}(s,I)_{\lambda}|_{c}=id_{T_{c}N}$ for all
$c\in M\cap U(C(I^{+}))_{1}$ and $r_{M}\circ\Psi^{M}(s,I)_{\lambda}=\psi
^{M}(s,I)_{\lambda}$. We define $\Phi(s,I)_{\lambda}:$ $TN|_{M}\rightarrow
TN|_{M}$ by $\Phi(s,I)_{\lambda}=(\Psi^{M}(s,I)_{\lambda})^{-1}$.

\section{Lemmas}

Let $I$ be the symbol in Theorem 4.1. In the proof of the following lemma,
$\Phi(s,I)_{\lambda}|_{c}$ ($c\in M$) is regarded as a linear isomorphism of
$T_{c}N$. Let $r_{0}$ be a small positive real number with $r_{0}<1/10$.

\begin{lemma}
Let $s\in\Gamma_{\Omega}^{tr}(N,P)$ be a section satisfying the hypotheses of
Theorem $4.1$. Then there exists a homotopy $s_{\lambda}$ relative to
$U(C(I^{+}))_{2-3r_{0}}$ in $\Gamma_{\Omega}^{tr}(N,P)$ with $s_{0}=s$ satisfying

$(1)$ for any $\lambda$, $S^{I}(s_{\lambda})=S^{I}(s)$ and $\pi_{P}^{\infty
}\circ s_{\lambda}|S^{I}(s_{\lambda})=\pi_{P}^{\infty}\circ s|S^{I}(s),$

$(2)$ for any point $c\in S^{I}(s_{1})$, we have $i_{K_{1}(S^{I}(s_{1}))}%
=\Psi(s_{1},I)$. In particular, $K_{1}(S^{I}(s_{1}))_{c}\subset\frak{n}%
(s,I)_{c}$.
\end{lemma}

\begin{proof}
Recall the exponential map $\exp_{N,x}:T_{x}N\rightarrow N$ defined near
$\mathbf{0}\in T_{x}N$. We write an element of $\frak{n}(I)_{c}$ as
$\mathbf{v}_{c}$. There exists a small positive number $\delta$ such that the
map
\[
e:D_{\delta}(\frak{n}(I))|_{M}\rightarrow N
\]
defined by $e(\mathbf{v}_{c})=\exp_{N,c}(\mathbf{v}_{c})$ is an embedding,
where $c\in M$ and $\mathbf{v}_{c}\in D_{\delta}(\frak{n}(I)_{c})$ (note that
$e|M$ is the inclusion). Let $\rho:[0,\infty)\rightarrow\mathbf{R}$ be a
decreasing smooth function such that $0\leq\rho(t)\leq1$, $\rho(t)=1$ if
$t\leq\delta/10$ and $\rho(t)=0$ if $t\geq\delta$.

If we represent $s(x)\in\Omega(N,P)$ by a jet $j_{x}^{\infty}\sigma_{x}$ for a
germ $\sigma_{x}:(N,x)\rightarrow(P,\pi_{P}^{\infty}\circ s(x))$, then we
define the homotopy $s_{\lambda}$ of $\Gamma_{\Omega}^{tr}(N,P)$ using
$\Phi(s,I)_{\lambda}$ by%
\begin{equation}
\left\{
\begin{array}
[c]{ll}%
\begin{array}
[c]{l}%
s_{\lambda}(e(\mathbf{v}_{c}))\\
=j_{e(\mathbf{v}_{c})}^{\infty}(\sigma_{e(\mathbf{v}_{c})}\circ\exp_{N,c}%
\circ\Phi(s,I)_{\rho(\Vert\mathbf{v}_{c}\Vert)\lambda}|_{c}\circ\exp
_{N,c}^{-1})
\end{array}
& \text{\textrm{if} $c\in M$ \textrm{and} }\Vert\text{$\mathbf{v}_{c}\Vert
\leq\delta,$}\\%
\begin{array}
[c]{l}%
s_{\lambda}(x)=s(x)
\end{array}
& \text{\textrm{if} $x\notin\mathrm{Im}(e).$}%
\end{array}
\right.
\end{equation}
Here, $\Phi(s,I)_{\rho(\Vert\mathbf{v}_{c}\Vert)\lambda}|_{c}$\ refers to
$\ell(\mathbf{v}_{c})\circ(\Phi(s,I)_{\rho(\Vert\mathbf{v}_{c}\Vert)\lambda
}|_{c})\circ\ell(-\mathbf{v}_{c})$, where $\ell(\mathbf{v})$ denotes the
parallel translation defined by $\ell(\mathbf{v})(\mathbf{a})=\mathbf{a}%
+\mathbf{v}$. If $\Vert\mathbf{v}_{c}\Vert\geq\delta$, then $\Phi
(s,I)_{\rho(\Vert\mathbf{v}_{c}\Vert)\lambda}|_{c}=\Phi(s,I)_{0}|_{c}$, and if
$c\in S^{I}\cap U(C(I^{+}))_{2-3r_{0}}$, then $\Phi(s,I)_{\lambda}|_{c}%
=\Phi(s,I)_{0}|_{c}$. Hence, $s_{\lambda}$ is well defined. It follows from
(5.1) that

(1) $\pi_{P}^{\infty}\circ s_{\lambda}(x)=\pi_{P}^{\infty}\circ s(x)$,

(2) $S^{I}(s_{\lambda})=S^{I}(s)$,

(3) if $c\in S^{I}(s)$, then we have that $\frak{n}(s,I)_{c}\supset
K_{1}(S^{I}(s_{1}))_{c}$ and $i_{K_{1}(S^{I}(s_{1}))}=\Psi(s_{1},I)$.

(4) $s_{\lambda}$ is transverse to $\Sigma^{I}(N,P)$.
\end{proof}

In what follows we set $d_{1}(s,I)=(s|S^{I})^{\ast}(\mathbf{d}_{1})$. We also
choose and fix a Riemannian metric of $P$ and identify $Q(S^{I}(s))$ with the
orthogonal complement of Im$(d_{1}(s,I))$ in $(\pi_{P}^{\infty}\circ
s|S^{I})^{\ast}(TP)$.

\begin{lemma}
Let $s$ be a section of $\Gamma_{\Omega}^{tr}(N,P)$ satisfying the property
$(2)$ for $s$ $($in place of $s_{1}$$)$ of Lemma $5.1$. Then there exists a
homotopy $s_{\lambda}$ relative to $U(C(I^{+}))_{2-3r_{0}}$ in $\Gamma
_{\Omega}^{tr}(N,P)$ with $s_{0}=s$ such that

$(1)$ $S^{I}(s_{\lambda})=S^{I}(s)$ for any $\lambda$,

$(2)$ $\pi_{P}^{\infty}\circ s_{1}|S^{I}(s_{1})$ is an immersion into $P$ such
that $d(\pi_{P}^{\infty}\circ s_{1}|S^{I}(s_{1})):T(S^{I}(s_{1}))\linebreak
\rightarrow TP$ is equal to $(\pi_{P}^{\infty}\circ s_{1})^{TP}\circ
d_{1}(s_{1},I)|T(S^{I}(s_{1}))$, where $(\pi_{P}^{\infty}\circ s_{1}%
)^{TP}:(\pi_{P}^{\infty}\circ s_{1})^{\ast}(TP)\rightarrow TP$ is the
canonical induced bundle map.

\begin{proof}
Since $\mathbf{K}_{1}\cap T(\Sigma^{I}(N,P))=\{\mathbf{0}\},$ it follows that
$(\pi_{P}^{\infty}\circ s)^{TP}\circ d_{1}(s,I)|T(S^{I})$ is a monomorphism.
By the Hirsch Immersion Theorem ([H1, Theorem 5.7]) there exists a homotopy of
monomorphisms $m_{\lambda}^{\prime}:T(S^{I})\rightarrow TP$ covering a
homotopy $m_{\lambda}:S^{I}\rightarrow P$ such that $m_{0}^{\prime}=(\pi
_{P}^{\infty}\circ s)^{TP}\circ d_{1}(s,I)|T(S^{I})$ and that $m_{1}$ is an
immersion with $d(m_{1})=m_{1}^{\prime}$. Then we can extend $m_{\lambda
}^{\prime}$ to a homotopy $\widetilde{m_{\lambda}^{\prime}}:TN|_{S^{I}%
}\rightarrow TP$ of homomorphisms of constant rank $n-i_{1}$ relative to
$U(C(I^{+}))_{2-3r_{0}}$. In fact, let $m:S^{I}\times\lbrack0,1]\rightarrow
P\times\lbrack0,1]$ and $m^{\prime}:T(S^{I})\times\lbrack0,1]\rightarrow
TP\times\lbrack0,1]$ be the maps defined by $m(c,\lambda)=(m_{\lambda
}(c),\lambda)$ and $m^{\prime}(\mathbf{v},\lambda)=(m_{\lambda}^{\prime
}(\mathbf{v}),\lambda)$ respectively. Let $m^{\ast}(m^{\prime}):T(S^{I}%
)\times\lbrack0,1]\rightarrow m^{\ast}(TP\times\lbrack0,1])$ be the canonical
monomorphism induced from $m^{\prime}$ by $m$. Let $\mathcal{F}_{1}%
=$Im$(m^{\ast}(m^{\prime}))$ and $\mathcal{F}_{2}$ be the orthogonal
complement of $\mathcal{F}_{1}$ in $m^{\ast}(TP\times\lbrack0,1])$. Since
$\mathcal{F}_{2}$ is isomorphic to $(\mathcal{F}_{2}|_{S^{I}\times0}%
)\times\lbrack0,1]$, we obtain a monomorphism of rank $c(I)-i_{1}$
\[
j_{\mathcal{F}}:\text{Im}(d_{1}(s,I)|\frak{n}(I))\times\lbrack0,1]\rightarrow
\mathcal{F}_{2}\text{ \ \ \ \ \ over }S^{I}\times\lbrack0,1]\text{.}%
\]
Since $d_{1}(s,I)|(TN|_{S^{I}})$ is of constant rank $n-i_{1}$, it induces the
homomorphism of kernel rank $i_{1}$%
\[
d:\frak{n}(I)\times\lbrack0,1]\rightarrow\text{Im}(d_{1}(s,I)|\frak{n}%
(I))\times\lbrack0,1]\overset{j_{\mathcal{F}}}{\rightarrow}\mathcal{F}%
_{2}\text{.}%
\]
We define $\widetilde{m^{\prime}}$ to be the composition%
\[%
\begin{array}
[c]{c}%
TN|_{S^{I}}\times\lbrack0,1]\cong(T(S^{I})\oplus\frak{n}(I))\times
\lbrack0,1]\overset{\underrightarrow{\text{ }m^{\ast}(m^{\prime})\oplus
d\text{ \ }}}{}\mathcal{F}_{1}\oplus\mathcal{F}_{2}\\
\rightarrow\text{Im}(m^{\ast}(m^{\prime}))\oplus\text{Cok}(m^{\ast}(m^{\prime
}))\cong m^{\ast}(TP\times\lbrack0,1])\overset{\underrightarrow{\text{
}m^{TP\times\lbrack0,1]}\text{\ }}}{}TP\times\lbrack0,1].
\end{array}
\]
We define $\widetilde{m_{\lambda}^{\prime}}$ to be $(\widetilde{m_{\lambda
}^{\prime}}(\mathbf{v}),\lambda)=\widetilde{m^{\prime}}(\mathbf{v},\lambda)$.

Next we construct a homotopy $s_{\lambda}:N\rightarrow\Omega(N,P)$ from
$\widetilde{m_{\lambda}^{\prime}}.$ Recall the submanifold $\widetilde{\Sigma
}^{i_{1}}(N,P)$ of $J^{1}(N,P)=J^{1}(TN,TP)$ which corresponds to
$\Sigma^{i_{1}}(N,P)$ in Section 2 (9). Then $\pi_{1}^{\infty}|\Sigma
^{I}(N,P):\Sigma^{I}(N,P)\rightarrow\widetilde{\Sigma}^{i_{1}}(N,P)$ becomes a
fiber bundle. We regard $\widetilde{m_{\lambda}^{\prime}}$ as a homotopy
$S^{I}\rightarrow\widetilde{\Sigma}^{i_{1}}(N,P).$\ By the covering homotopy
property to $s|S^{I}$and $\widetilde{m_{\lambda}^{\prime}}$, we obtain a
homotopy $s_{\lambda}^{\prime}:S^{I}\rightarrow\Sigma^{I}(N,P)$ covering
$\widetilde{m_{\lambda}^{\prime}}$ relative to $U(C(I^{+}))_{2-3r_{0}}$ such
that $s_{0}^{\prime}=s|S^{I}$.

By using the transversality of $s$ and the homotopy extension property to $s$
and $s_{\lambda}^{\prime}$, we first extend $s_{\lambda}^{\prime}$ to a
homotopy defined on a tubular neighborhood of $S^{I}$ and then to a required
homotopy $s_{\lambda}\in\Gamma_{\Omega}^{tr}(N,P)$, which satisfies $s_{0}=s$,
$s_{\lambda}|S^{I}=s_{\lambda}^{\prime}$ and $\ s_{\lambda}|U(C(I^{+}%
))_{2-3r_{0}}=s|U(C(I^{+}))_{2-3r_{0}}$.
\end{proof}
\end{lemma}

Here we give two lemmas necessary for the proof of Theorem 4.1. Let
$\pi:E\rightarrow S$ be a smooth $c(I)$-dimensional vector bundle with a
metric over an $(n-c(I))$-dimensional manifold, where $S$ is identified with
the zero-section. Then we can identify $\exp_{E}|D_{\varepsilon}%
(E)=id_{D_{\varepsilon}(E)}$.

\begin{lemma}
Let $\pi:E\rightarrow S$ be given as above. Let $f_{i}:E\rightarrow P$
$(i=1,2)$ be $\Omega$-regular maps which have singularities of the symbol $I$
exactly on $S$\ such that

$(\mathrm{i})$ $f_{1}|S=f_{2}|S$, which are immersions,

$(\mathrm{ii})$ $S=S^{I}(j^{\infty}f_{1})=S^{I}(j^{\infty}f_{2})$,

$(\mathrm{iii})$ $K_{1}(S^{I}(j^{\infty}f_{1}))_{c}=K_{1}(S^{I}(j^{\infty
}f_{2}))_{c}$ are tangent to $\pi^{-1}(c)$,

$(\mathrm{iv})$ $T_{c}(S^{I_{j-1}}(j^{\infty}f_{1}))=T_{c}(S^{I_{j-1}%
}(j^{\infty}f_{2}))$, $((j^{\infty}f_{1})^{\ast}\mathbf{P}_{j})_{c}%
=((j^{\infty}f_{2})^{\ast}\mathbf{P}_{j})_{c}$ and $((j^{\infty}f_{1})^{\ast
}(\mathbf{d}_{j+1}\circ d(j^{\infty}f_{1})))_{c}=((j^{\infty}f_{2})^{\ast
}(\mathbf{d}_{j+1}\circ d(j^{\infty}f_{2})))_{c}$ for each number $j$ and any
$c\in S$.

\noindent Let $\eta:S\rightarrow\lbrack0,1]$ be any smooth function. Let
$\varepsilon:S\rightarrow\mathbf{R}$ be a sufficiently small positive smooth
function. Let $\mathbf{f}^{\eta}(\mathbf{v}_{c})$ denote $\exp_{P,f_{1}%
(c)}((1-\eta(c))\exp_{P,f_{1}(c)}^{-1}(f_{1}(\mathbf{v}_{c}))+\eta
(c)\exp_{P,f_{2}(c)}^{-1}(f_{2}(\mathbf{v}_{c})))$ for any $c\in S$ and any
$\mathbf{v}_{c}\in\pi^{-1}(c)$ with $\Vert\mathbf{v}_{c}\Vert\leq
\varepsilon(c)$.

Then the map $\mathbf{f}^{\eta}:D_{\varepsilon}(E)\rightarrow P$ is a
well-defined $\Omega$-regular map such that

$(1)$ $\mathbf{f}^{\eta}|S=f_{1}|S=f_{2}|S,$

$(2)$ $S=S^{I}(j^{\infty}\mathbf{f}^{\eta})$,

$(3)$ $T_{c}(S^{I_{j-1}}(j^{\infty}\mathbf{f}^{\eta}))=T_{c}(S^{I_{j-1}%
}(j^{\infty}f_{1})),$ $((j^{\infty}\mathbf{f}^{\eta})^{\ast}\mathbf{P}%
_{j})_{c}=((j^{\infty}f_{1})^{\ast}\mathbf{P}_{j})_{c}$ and $((j^{\infty
}\mathbf{f}^{\eta})^{\ast}(\mathbf{d}_{j+1}\circ d(j^{\infty}\mathbf{f}^{\eta
})))_{c}=((j^{\infty}f_{1})^{\ast}(\mathbf{d}_{j+1}\circ d(j^{\infty}%
f_{1})))_{c}$ for each number $j$ and any $c\in S$,

\begin{proof}
Let us take a Riemannian metric on $E$ which is compatible with the metric of
the vector bundle $E$\ over $S$. In particular, $S$\ is a Riemannian
submanifold\ of $E$. Furthermore, take a Riemannian metric on $P$\ such that
$f_{i}(S)\cap P$\ is a Riemannian submanifold\ of $P$\ around $f(c)$. Then the
local coordinates of $\exp_{N,c}(K_{1,c})$\ and $\exp_{P,f_{i}(c)}(Q_{c}%
)$\ are independent of the coordinates of $S$, where $Q_{c}$\ is regarded as a
subspace of $T_{f(c)}P$.\ 

We may consider $\eta(c)$ as a constant when dealing with higher intrinsic
derivatives in the lemma by the identification (1.2) and the property of the
total tangent bundle $\mathbf{D}$ given in the beginning\ of Section 2. Then
the assertion follows from the assumptions and the properties of
$\Sigma^{I_{j}}(N,P)$.
\end{proof}
\end{lemma}

The proof of the following lemma is elementary, and so is left to the reader.

\begin{lemma}
Let $\pi:E\rightarrow S$ be given as above. Let $(\Omega,\Sigma)$ be a pair of
a smooth manifold and its submanifold of codimension $c(I)$. Let
$\varepsilon:S\rightarrow\mathbf{R}$ be a sufficiently small positive smooth
function. Let $h:D_{\varepsilon}(E)\rightarrow(\Omega,\Sigma)$ be a smooth map
such that $S=h^{-1}(\Sigma)$ and that $h$ is transverse to $\Sigma$. Then
there exists a smooth homotopy $h_{\lambda}:(D_{\varepsilon}(E),S)\rightarrow
(\Omega,\Sigma)$ between $h$ and $\exp_{\Omega}\circ dh|D_{\varepsilon}(E)$
such that

$(1)$ $h_{\lambda}|S=h_{0}|S$, $S=h_{\lambda}^{-1}(\Sigma)=h_{0}^{-1}(\Sigma)$
for any $\lambda$,

$(2)$ $h_{\lambda}$ is smooth and is transverse to $\Sigma$ for any $\lambda$,

$(3)$ $h_{0}=h$ and $h_{1}(\mathbf{v}_{c})=\exp_{\Omega,h(c)}\circ
dh(\mathbf{v}_{c})$ for $c\in S$ and $\mathbf{v}_{c}\in D_{\varepsilon}(E_{c}).$
\end{lemma}

\section{Proof of Theorem 4.1}

Consider the bundles $\frak{n}(I)$ ($=\frak{n}(s,I)$) and $Q$ $(=Q(S^{I}(s))$.
For a point $c\in S^{I}(s)$, take an open neighborhood $U$ around $c$ such
that $\frak{n}(I)|_{U}$ and $Q|_{U}$ are the trivial bundles, say
$U\times\mathbf{R}^{c(I)}$ and $U\times\mathbf{R}^{p-n+i_{1}}$ respectively,
where $\mathbf{R}^{c(I)}$ has coordinates $(x_{1},\ldots,x_{c(I)})$ and
$\mathbf{R}^{p-n+i_{1}}$ has $(y_{1},\ldots,y_{p-n+i_{1}})$. Then we identify
an element of $\mathrm{Hom}(\bigcirc^{j}\frak{n}(I),Q)|_{U}$ with polynomials
$y_{i}(c)=\Sigma_{|\omega|=j}a_{i}^{\omega}(c)x_{1}^{\omega_{1}}x_{2}%
^{\omega_{2}}\cdots x_{c(I)}^{\omega_{c(I)}}$, $c\in U$, where $\omega
=(\omega_{1},\omega_{2},\cdots,\omega_{c(I)})$, $\omega_{\ell}\geq0$
($i=1,\cdots,p-n+i_{1}$), and $|\omega|=\omega_{1}+\cdots+\omega_{c(I)}$ and
$a_{i}^{\omega}(c)$ are real numbers. If $a_{i}^{\omega}(c)$ are smooth
functions of $c,$ then $\{a_{i}^{\omega}(c)\}$ defines a smooth section of
Hom$(\bigcirc^{j}\frak{n}(I),Q)|_{S^{I}}$ over $U$.

We first introduce several homomorphisms between vector bundles over
$S^{I}(s)$, which are used for the construction of the required $\Omega
$-regular map in Theorem 4.1.

Let $s\in\Gamma_{\Omega}^{tr}(N,P)$. By deforming $s$ if necessary, we may
assume without loss of generality that $s$ satisfies (2) of Lemma 5.1 and (2)
of Lemma 5.2, where $s_{1}$ is replaced by $s$.

In the following, let $K_{j}$ ($j\geq1$) refer to $K_{j}(S^{I}(s))$. Let
$K_{j}/K_{j+1}$ refer to the orthogonal complement of $K_{j+1}$ in $K_{j}$,
$T_{j}^{I}$ refer to the orthogonal complement of $K_{j}/K_{j+1}$ in
$\frak{n}(s,I_{j}\subset I_{j-1})$, and $P_{j}^{I}=(s|S^{I})^{\ast}%
\mathbf{P}_{j}$. Then we have the following isomorphism by (2.2)%
\begin{align}
(s|S^{I})^{\ast}(\mathbf{d}_{j+1}\circ ds|\frak{n}(s,I_{j}  &  \subset
I_{j-1})):\\
\frak{n}(s,I_{j}  &  \subset I_{j-1})=K_{j}/K_{j+1}\oplus T_{j}^{I}\rightarrow
P_{j}^{I}\text{ }(1\leq j\leq k).\nonumber
\end{align}

We first define the section $\frak{q}^{\prime}(s,I)^{1}$ of $\mathrm{Hom}%
(\frak{n}(I),\mathrm{\operatorname{Im}}(d_{1}(s,I)))$ over $S^{I}(s)$ defined
by $\frak{q}^{\prime}(s,I)^{1}=d_{1}(s,I)|\frak{n}(I)$, which vanishes on
$K_{1}|_{S^{I}}$ and gives an isomorphism of $\oplus_{j=1}^{k}T_{j}^{I}$ onto
$\mathrm{\operatorname{Im}}(d_{1}(s,I))$.

For $1\leq j\leq k$, $\mathbf{q}(k)^{j+1,j+1}$ in (3.1) induces the
homomorphism%
\begin{equation}
\frak{q}(s,I)^{j+1}:\frak{n}(s,I_{j}\subset I_{j-1})\bigcirc K_{j}\bigcirc
K_{j-1}\bigcirc\cdots\bigcirc K_{1}\rightarrow Q
\end{equation}
over $S^{I}(s)$ defined as the composition%
\[
((s|S^{I})^{\ast}\mathbf{q}(k)^{j+1,j+1})\circ(((s|S^{I})^{\ast}%
ds|\frak{n}(s,I_{j}\subset I_{j-1}))\bigcirc id_{K_{j}\bigcirc K_{j-1}%
\bigcirc\cdots\bigcirc K_{1}}).
\]
Furthermore, we define the following section of Hom$(\Sigma_{j=1}%
^{k+1}\bigcirc^{j}\frak{n}(I),Q)$%
\begin{equation}%
\begin{array}
[c]{ll}%
\frak{q}^{\prime}(s,I)=\Sigma_{j=1}^{k}\frak{q}(s,I)^{j+1} & \text{over }%
S^{I}(s).
\end{array}
\end{equation}
Let us consider the direct sum decompositions%
\begin{align*}
\frak{n}(s,I)  &  =\oplus_{j=1}^{k}\frak{n}(s,I_{j}\subset I_{j-1}),\text{
\ \ }\frak{n}(s,I_{j}\subset I_{j-1})=K_{j}/K_{j+1}\oplus T_{j}^{I},\\
K_{1}  &  =\oplus_{j=1}^{k-1}K_{j}/K_{j+1}\oplus K_{k},\text{ \ \ }(\pi
_{P}^{\infty}\circ s|S^{I})^{\ast}(TP)=Q\oplus Q^{\bot}%
\end{align*}
and the inclusion $i_{Q}:Q\rightarrow(\pi_{P}^{\infty}\circ s|S^{I})^{\ast
}(TP)$. Then we obtain the smooth\ fiber map
\begin{equation}
\frak{q}(s,I)=(\pi_{P}^{\infty}\circ s|S^{I})^{TP}\circ(i_{Q}\circ
\frak{q}^{\prime}(s,I)+\frak{q}^{\prime}(s,I)^{1}):\frak{n}(s,I)\rightarrow TP
\end{equation}
covering the immersion $\pi_{P}^{\infty}\circ s|S^{I}(s):S^{I}(s)\rightarrow
P$ such that for any $c\in S^{I}(s)$, $\frak{q}(s,I)_{c}$ is regarded as
$p-n+c(I)$ polynomials of $c(I)$ variables with constant $0$.

\begin{proof}
[Proof of Theorem 4.1]By Lemmas 5.1 and 5.2 we may assume that $s$ satisfies
(2) of Lemma 5.1 and (2) of Lemma 5.2, where $s_{1}$ is replaced by $s$. We
define $E(S^{I})$ to be the union of all $\exp_{N}(D_{\delta\circ s}%
(\frak{n}(I)))$, where $\delta:\Sigma^{I}(N,P)\rightarrow\mathbf{R}$ is a
sufficiently small positive function such that $\delta\circ s|(S^{I}\setminus
$Int$U(C(I^{+})){_{2}})$ is constant. This is a tubular neighborhood of
$S^{I}$.

It is enough for the proof of Theorem 4.1 to prove the following assertion:

(\textbf{A}) {There exists a homotopy $H_{\lambda}$ relative to $U(C(I^{+}%
))_{2-r_{0}}$ in $\Gamma_{\Omega}^{tr}(N,P)$ with $H_{0}=s$ satisfying the
following. }

{$(1)$ $S{^{I}}(H_{\lambda})=S^{I}$ for} {any $\lambda$. }

{$(2)$ We have an }$\Omega$-regular{ map $G$ defined on a neighborhood of
${U(}${$C(I^{+})$}${)}_{2-r_{0}}$ to $P$ such that $j^{\infty}G=H_{1}$ on
}${U(}${$C(I^{+})$}${)}_{2-r_{0}}\cup E(S^{I}).$

By the Riemannian metric on $P$, we identify $Q$ with the orthogonal
$p-n+i_{1}$ dimensional bundle of Im$(d_{1}(s,I))$ in $(\pi_{P}^{\infty}\circ
s|S^{I})^{\ast}(TP)$. Then the map $\exp_{P}\circ(\pi_{P}^{\infty}\circ
s|S^{I})^{TP}|D_{\gamma}(Q)$ is an immersion for some small positive function
$\gamma$. In the proof we express a point of $E(S^{I})$ as $\mathbf{v}_{c}$,
where $c\in S^{I},$ $\mathbf{v}_{c}\in\frak{n}(I)_{c}$ and $\Vert
\mathbf{v}_{c}\Vert\leq\delta(s(c))$. In the proof we say that a smooth
homotopy
\[
k_{\lambda}:(E(S^{I}),\partial E(S^{I}))\rightarrow(\Omega(N,P),\Omega
(N,P)\setminus\Sigma^{I}(N,P))
\]
has the property (C) if it satisfies that for any $\lambda$

(C-1) $k_{\lambda}^{-1}(\Sigma^{I}(N,P))=S^{I}$, and $\pi_{P}^{\infty}\circ
k_{\lambda}|S^{I}=\pi_{P}^{\infty}\circ k_{0}|S^{I\text{ }}$and,

(C-2) $k_{\lambda}$ is smooth and transverse to $\Sigma^{I}(N,P)$.

If we choose $\delta$ sufficiently small compared with $\gamma$, then we can
define the $\Omega$-regular map $g_{0}:E(S^{I})\rightarrow P$ by
\begin{equation}
g_{0}(\mathbf{v}_{c})=\exp_{P,\pi_{P}^{\infty}\circ s(c)}\circ\frak{q}%
(s,I)_{c}\circ\exp_{N,c}^{-1}(\mathbf{v}_{c}).
\end{equation}
It follows from Section 2 that $g_{0}$ has each point $c\in S^{I}$ as a
singularity of the symbol $I$\ and vice versa. Now we need to modify $g_{0}$
by using Lemma 5.3 so that $g_{0}$ is compatible with $g(I^{+})$. Let
$\eta:S^{I}\rightarrow\mathbf{R}$ be a smooth function such that

(i) $0\leq\eta(c)\leq1$ for $c\in S^{I},$

(ii) $\eta(c)=0$ for $c\in S^{I}\cap{U(}${$C(I^{+})$}${)}_{2-3r_{0}}$,

(iii) $\eta(c)=1$ for $c\in S^{I}\setminus{U(}${$C(I^{+})$}${)}_{2-4r_{0}}$.

\noindent Then consider the map $G:{U(}${$C(I^{+})$}${)}_{2-3r_{0}}\cup
E(S^{I})\rightarrow P$ defined by
\[
\left\{
\begin{array}
[c]{ll}%
G(x)=g(I^{+})(x) & \text{if $x\in{U(}${$C(I^{+})$}${)}_{2-3r_{0}}$},\\
G(\mathbf{v}_{c})=(1-\eta(c))g(I^{+})(\mathbf{v}_{c})+\eta(c)g_{0}%
(\mathbf{v}_{c}) & \text{if }\mathbf{v}_{c}\text{$\in E(S^{I})$}.
\end{array}
\right.
\]
It follows from Lemma 5.3 that $G$ is an $\Omega$-regular map defined ${U(}%
${$C(I^{+})$}${)}_{2-3r_{0}}\cup E(S^{I})$, that $G|E(S^{I})$ has the
singularities of the symbol $I$ exactly on $S^{I}$, and that for any $c\in
S^{I}$, the assumptions (i)-(iv) of Lemma 5.3 holds for $f_{1}=g(I^{+})$ and
$f_{2}=g_{0}$.

Set $\exp_{\Omega}=\exp_{\Omega(N,P)}$ for short. Let $h_{1}^{1}$ and
$h_{0}^{3}$ be the maps $(E(S^{I}),S^{I})\rightarrow(\Omega(N,P),\Sigma
^{I}(N,P))$ defined by
\begin{align*}
h_{1}^{1}(\mathbf{v}_{c})  &  =\exp_{\Omega,s(c)}\circ d_{c}s\circ(\exp
_{N,c})^{-1}(\mathbf{v}_{c}),\\
h_{0}^{3}(\mathbf{v}_{c})  &  =\exp_{\Omega,j^{\infty}G(c)}\circ
d_{c}(j^{\infty}G)\circ(\exp_{N,c})^{-1}(\mathbf{v}_{c}).
\end{align*}
By applying Lemma 5.4 to the section $s$ and $h_{1}^{1}$, we first obtain a
homotopy $h_{\lambda}^{1}\in\Gamma_{\Omega}^{tr}(E(S^{I}),P)$ between
$h_{0}^{1}=s$ and $h_{1}^{1}$ on $E(S^{I})$ satisfying the properties (1), (2)
and (3) of Lemma 5.4. Similarly we obtain a homotopy $h_{\lambda}^{3}\in
\Gamma_{\Omega}^{tr}(E(S^{I}),P)$ between $h_{0}^{3}$ and $h_{1}^{3}%
=j^{\infty}G$ on $E(S^{I})$ satisfying the properties (1), (2) and (3) of
Lemma 5.4.

Next we construct a homotopy of bundle maps $\frak{n}(I)\rightarrow\nu
(\Sigma^{I})$ covering a homotopy $S^{I}\rightarrow\Sigma^{I}(N,P)$ between
$ds|\frak{n}(I)$ and $d(j^{\infty}G)|\frak{n}(I)$. Let us recall the additive
structure of $J^{\infty}(N,P)$ in (1.2). Then we have the homotopy
$\kappa_{\lambda}:S^{I}\rightarrow J^{\infty}(N,P)$ defined by
\[
\kappa_{\lambda}(c)=(1-\lambda)s(c)+\lambda j^{\infty}G(c)\quad\text{covering
}\pi_{P}^{\infty}\circ s|S^{I}:S^{I}\rightarrow P,
\]
where $\pi_{P}^{\infty}\circ s|S^{I}$ is the immersion as in $(2)$ of Lemma
5.2. We show that $\kappa_{\lambda}$ is actually a homotopy of $S^{I}$ into
$\Sigma^{I}(N,P)$. Under the identification $(s)^{\ast}\mathbf{P}\cong(\pi
_{P}^{\infty}\circ s)^{\ast}(TP)$ and $s^{\ast}\mathbf{D}\cong TN$, it follows
from the decomposition of $\frak{n}(I)$ in (4.3) that%
\begin{equation}
(s|S^{I})^{\ast}(\mathbf{d}_{j+1}\circ ds|\frak{n}(I_{j}\subset I_{j-1}%
))=(j^{\infty}G|S^{I})^{\ast}(\mathbf{d}_{j+1}\circ d(j^{\infty}%
G)|\frak{n}(I_{j}\subset I_{j-1}))
\end{equation}
over $S^{I}$. These formulas are the direct consequence of the construction of
$\frak{q}(s,I)$ used in the definition of $G$ and the definition of the
intrinsic derivatives in Sections 2 and 3. By (6.6) we have that
$\frak{n}(\kappa_{\lambda},I)_{c}=\frak{n}(I)_{c}$ and $Q(\kappa_{\lambda
})_{c}=Q_{c}$ for any $c\in S^{I}$. Hence, it follows that the equalities of
the homomorphisms in (6.6) also hold when $s$ is replaced by $\kappa_{\lambda
}$ $(0\leq\lambda\leq1)$. This implies that $\kappa_{\lambda}$\ is a homotopy
into $\Sigma^{I}(N,P)$. Hence, the homotopy $(\kappa_{\lambda})^{\nu
(\Sigma^{I})}:\kappa_{\lambda}^{\ast}(\nu(\Sigma^{I}))\rightarrow\nu
(\Sigma^{I})$, $ds$ and $d(j^{\infty}G)$\ induce the homotopy of bundle maps
$\widetilde{\kappa_{\lambda}}:\frak{n}(I)\rightarrow\nu(\Sigma^{I})$ covering
$\kappa_{\lambda}$\ such that $\widetilde{\kappa_{0}}=ds$ and $\widetilde
{\kappa_{1}}=d(j^{\infty}G)$. We define the homotopy $h_{\lambda}^{2}%
:(E(S^{I}),S^{I})\rightarrow(\Omega(N,P),\Sigma^{I}(N,P))$ by%
\[
h_{\lambda}^{2}(\mathbf{v}_{c})=\exp_{\Omega,s(c)}\circ\widetilde
{\kappa_{\lambda}}\circ(\exp_{N,c})^{-1}(\mathbf{v}_{c}).
\]
Then we have that $h_{0}^{2}(\mathbf{v}_{c})=h_{1}^{1}(\mathbf{v}_{c}%
)=\exp_{\Omega,s(c)}\circ d_{c}s\circ(\exp_{N,c})^{-1}(\mathbf{v}_{c})$ and
$h_{0}^{3}(\mathbf{v}_{c})=h_{1}^{2}(\mathbf{v}_{c})=\exp_{\Omega,j^{\infty
}G(c)}\circ d_{c}(j^{\infty}G)\circ(\exp_{N,c})^{-1}(\mathbf{v}_{c})$ on
$E(S^{I})$.

Since $h_{0}^{1}(\mathbf{v}_{c})=h_{1}^{3}(\mathbf{v}_{c})=s(\mathbf{v}_{c})$
for $\mathbf{v}_{c}\in E(S^{I})\cup{U(}${$C(I^{+})$}${)}_{2-3r_{0}}$, we may
assume in the construction of $h_{\lambda}^{1}$, $h_{\lambda}^{2}$ and
$h_{\lambda}^{3}$ that if $\mathbf{v}_{c}\in E(S^{I})\cup{U(}${$C(I^{+})$}%
${)}_{2-3r_{0}}$, then $h_{\lambda}^{2}(\mathbf{v}_{c})=h_{0}^{2}%
(\mathbf{v}_{c})=h_{1}^{2}(\mathbf{v}_{c})$ and $h_{\lambda}^{1}%
(\mathbf{v}_{c})=h_{1-\lambda}^{3}(\mathbf{v}_{c})$ for any $\lambda$.

Let $\overline{h}_{\lambda}\in\Gamma_{\Omega}^{tr}(E(S^{I})\cup{U(}$%
{$C(I^{+})$}${)}_{2-3r_{0}},P)$ be the homotopy which is obtained by pasting
$h_{\lambda}^{1}$, $h_{\lambda}^{2}$ and $h_{\lambda}^{3}$. The homotopies
$h_{\lambda}^{1}$ and $h_{\lambda}^{3}$ are not homotopies relative to
$E(S^{I})\cap{U(}${$C(I^{+})$}${)}_{2-3r_{0}}$ in general. By using the above
properties of $h_{\lambda}^{1}$, $h_{\lambda}^{2}$ and $h_{\lambda}^{3}$, we
can modify $\overline{h}_{\lambda}$ to a homotopy $h_{\lambda}\in
\Gamma_{\Omega}^{tr}(E(S^{I}),P)$\ satisfying the property (C) such that

(1) $h_{\lambda}(\mathbf{v}_{c})=h_{0}(\mathbf{v}_{c})=s(\mathbf{v}_{c})$ for
any $\lambda$\ and any $\mathbf{v}_{c}\in E(S^{I})\cap{U(}${$C(I^{+})$}%
${)}_{2-2r_{0}}$,

(2) $h_{0}(\mathbf{v}_{c})=s(\mathbf{v}_{c})$ for any $\mathbf{v}_{c}\in
E(S^{I}),$

(3) $h_{1}(\mathbf{v}_{c})=j^{\infty}G(\mathbf{v}_{c})$ for any $\mathbf{v}%
_{c}\in E(S^{I})$.

By (1), we can extend $h_{\lambda}$ to the homotopy $H_{\lambda}^{\prime}%
\in\Gamma_{\Omega}^{tr}(E(S^{I})\cup{U(}${$C(I^{+})$}${)}_{2-2r_{0}}%
,P)$\ defined by $H_{\lambda}^{\prime}|E(S^{I})=h_{\lambda}$ and $H_{\lambda
}^{\prime}|{U(}${$C(I^{+})$}${)}_{2-2r_{0}}=s|{U(}${$C(I^{+})$}${)}_{2-2r_{0}%
}$.

By the transversality of $H_{\lambda}^{\prime}$ and the homotopy extension
property to $s$ and $H_{\lambda}^{\prime}$, we obtain an extended homotopy
\[
H_{\lambda}:(N,S^{I})\rightarrow(\Omega(N,P),\Sigma^{I}(N,P))
\]
relative to ${U(}${$C(I^{+})$}${)}_{2-r_{0}}$\ with $H_{0}=s$. Furthermore, we
replace $\delta$ and $E(S^{I})$\ by smaller ones. Then $H_{\lambda}$ is a
required homotopy in $\Gamma_{\Omega}^{tr}(N,P)$ in the assertion (\textbf{A}).
\end{proof}

\section{Proof of Theorem 0.2}

In this section we prove Theorem 0.2 by applying Theorem 0.1.

\begin{proposition}
Under the same assumption of Theorem 0.2, any section $s\in\Gamma
_{\Omega^{I_{r}}}^{tr}(N,P)$ has a homotopy $s_{\lambda}\in\Gamma
_{\Omega^{I_{r}}}^{tr}(N,P)$ such that

$(1)$ $s_{0}=s,$

$(2)$ $s_{1}$ is a section of $\Omega^{J}(N,P)$ over $N$,

$(3)$ $S^{I_{r}}(s_{\lambda})=S^{I_{r}}(s)=S^{I_{r}}(s_{1})$ for any $\lambda$.
\end{proposition}

We need the following lemma for the proof of Proposition 7.1.

\begin{lemma}
Assume the same assumption of Theorem 0.2. Then, we have

$(1)$ $\mathbf{Q}_{1}\mathbf{|}_{\Sigma^{I_{r}}(N,P)}$, and $\bigcirc
^{2}\mathbf{K}_{r}|_{\Sigma^{I_{r}}(N,P)}$ are trivial line bundles equipped
with the canonical orientations respectively,

$(2)$ The homomorphisms $\mathbf{c}_{j}|\mathrm{Hom}(\bigcirc^{j}%
\mathbf{K}_{r},\mathbf{Q}_{1}):\mathrm{Hom}(\bigcirc^{j}\mathbf{K}%
_{r},\mathbf{Q}_{1})\rightarrow\mathbf{P}_{j}$\ $(1\leq j\leq r)$ and
$\mathbf{e}_{j-1}\circ\mathbf{c}_{j-1}:\mathrm{Hom}(\bigcirc^{j-1}%
\mathbf{K}_{r},\mathbf{Q}_{1})\rightarrow\mathbf{Q}_{j}$\ $(1<j\leq r)$ are
injective over $\Sigma^{I_{r}}(N,P)$.
\end{lemma}

\begin{proof}
(1) By Section 2 (5), $\mathbf{d}_{2}|\mathbf{K}_{1}:\mathbf{K}_{1}%
\rightarrow\mathrm{Hom}(\mathbf{K}_{1},\mathbf{Q}_{1})$ induces the
isomorphism%
\[
\mathbf{K}_{1}\mathbf{/K}_{2}\rightarrow\text{Hom}(\mathbf{K}_{1}%
\mathbf{/K}_{2},\mathbf{Q}_{1})\text{ \ \ \ over }\Sigma^{I_{r}}(N,P).
\]
This yields $\mathbf{q}:\mathbf{K}_{1}/\mathbf{K}_{2}\bigcirc\mathbf{K}%
_{1}/\mathbf{K}_{2}\rightarrow\mathbf{Q}_{1}$\ over $\Sigma^{I_{r}}(N,P)$,
which is a nonsingular quadratic form on each fiber. Since dim$\mathbf{K}%
_{1}\mathbf{/K}_{2}=n-p+1-i_{2}$ is odd, we choose the unique orientation of
$\mathbf{Q}_{1}$, expressed by the unit vector $\mathbf{e}_{p}$, so that the
index (the number of the negative eigen values) of $\mathbf{q}_{z}$,
$z\in\Sigma^{I_{r}}(N,P)$ is less than $(n-p+1-i_{2})/2$.

Since $\mathbf{K}_{r}|_{\Sigma^{I_{r}}(N,P)}$ is a line bundle, $\bigcirc
^{2}\mathbf{K}_{r}|_{\Sigma^{I_{r}}(N,P)}$ has the canonical orientation.

(2) We prove the assertion by induction on $j$ ($r\geq3$). Let $j=1$. Since
the kernel of $\mathbf{d}_{2}|\mathbf{K}_{1}$ is $\mathbf{K}_{2}$, we have
that $\mathbf{c}_{1}=\mathbf{u}_{1}$ induces the inclusion Hom$(\mathbf{K}%
_{r},\mathbf{Q}_{1})|_{\Sigma^{I_{r}}}\subset\mathrm{Hom}(\mathbf{K}%
_{1},\mathbf{Q}_{1})|_{\Sigma^{I_{r}}}=\mathbf{P}_{1}|_{\Sigma^{I_{r}}}$ and
$\mathbf{e}_{1}|_{\Sigma^{I_{r}}}:\mathbf{P}_{1}|_{\Sigma^{I_{r}}}%
\rightarrow\mathbf{Q}_{2}|_{\Sigma^{I_{r}}}$ is identified with the
restriction Hom$(\mathbf{K}_{1},\mathbf{Q}_{1})|_{\Sigma^{I_{r}}}%
\rightarrow\mathrm{Hom}(\mathbf{K}_{2},\mathbf{Q}_{1})|_{\Sigma^{I_{r}}}$.
Hence, $\mathbf{e}_{1}\circ\mathbf{c}_{1}|($Hom$(\mathbf{K}_{r},\mathbf{Q}%
_{1})|_{\Sigma^{I_{r}}})$ is injective. Suppose that $\mathbf{e}_{j-2}%
\circ\mathbf{c}_{j-2}|\mathrm{Hom}(\bigcirc^{j-2}\mathbf{K}_{r},\mathbf{Q}%
_{1})$ is injective into $\mathbf{Q}_{j-1}$ over $\Sigma^{I_{r}}(N,P)$ for
$j-2<r$. Then it follows from the definition of $\mathbf{u}_{j-1}$ that
$\mathbf{c}_{j-1}|\mathrm{Hom}(\bigcirc^{j-1}\mathbf{K}_{r},\mathbf{Q}_{1})$
is injective into $\mathrm{Hom}(\mathbf{K}_{j-1},\mathbf{Q}_{j-1})$ over
$\Sigma^{I_{r}}(N,P)$. Since the image of $\mathbf{c}_{j-1}$ is $\mathbf{P}%
_{j-1}$, the map $\mathbf{c}_{j-1}|\mathrm{Hom}(\bigcirc^{j-1}\mathbf{K}%
_{r},\mathbf{Q}_{1})$ is injective into $\mathbf{P}_{j-1}$. Since
$\mathbf{d}_{j}|\mathbf{K}_{r}$ vanishes for $j\leq r$ over $\Sigma^{I_{r}%
}(N,P)$ and since $\mathbf{d}_{j}|\mathbf{K}_{r}$ is symmetric, the
composition $\mathbf{e}_{j-1}\circ\mathbf{c}_{j-1}|\mathrm{Hom}(\bigcirc
^{j-1}\mathbf{K}_{r},\mathbf{Q}_{1})$ is injective into $\mathbf{Q}_{j}$ over
$\Sigma^{I_{r}}(N,P)$ for $j\leq r$. Thus the map $\mathbf{c}_{j}%
|\mathrm{Hom}(\bigcirc^{j}\mathbf{K}_{r},\mathbf{Q}_{1})$ is injective into
$\mathbf{P}_{j}$ for $j\leq r$. This proves the lemma.
\end{proof}

\begin{proof}
[Proof of Proposition 7.1]In the proof we identify $J^{k}(N,P)$ with
$J^{k}(TN,TP)$ by (1.2). By (9) in Section 2, there exists the open subbundles
$\widetilde{\Omega}^{L}(N,P)$ of $J^{k}(N,P)$ such that $(\pi_{k}^{\infty
})^{-1}(\widetilde{\Omega}^{L}(N,P))=\Omega^{L}(N,P)$ for $L$ with length $k$.
It follows that $(\pi_{r}^{\infty}\circ s)(N\setminus(S^{I_{r}}(s)))\subset
\widetilde{\Omega}^{I_{r-1},0}(N\setminus(S^{I_{r}}(s)),P)$.

We now construct a new section $\widetilde{u}:N\rightarrow\widetilde{\Omega
}^{J}(N,P)$ as follows.

Let $\mathbf{e}_{p}(Q)_{c}$ and $\mathbf{e}(\bigcirc^{r+1}(K_{r})_{c})$ be the
oriented vectors induced from $\mathbf{e}_{p}(\mathbf{Q}_{s(c)})$,
$\mathbf{e}(\bigcirc^{r+1}\mathbf{K}_{r,s(c)})$ by $s$ respectively. Then we
define the section $\phi^{J}:S^{I_{r}}(s)\rightarrow\mathrm{Hom}%
(\bigcirc^{r+1}K_{r},Q)$ by $\phi^{J}(c)(\mathbf{e}(\bigcirc^{r+1}%
(K_{r})_{s(c)}))=\mathbf{e}_{p}(Q_{s(c)})$. Then we can extend $\phi^{J}$\ to
a section $u_{\phi}:S^{I_{r}}(s)\rightarrow\mathrm{Hom}(S^{r+1}((\pi
_{N}^{\infty}\circ s)^{\ast}(TN)),(\pi_{P}^{\infty}\circ s)^{\ast}(TP))$ such
that $u_{\phi}(c)|\bigcirc^{r+1}K_{r,c}=\phi^{J}(c)$ for $c\in S^{I_{r}}(s)$.
Since $S^{I_{r}}(s)$ is a closed submanifold and since Hom$(S^{r+1}((\pi
_{N}^{\infty}\circ s)^{\ast}(TN)),(\pi_{P}^{\infty}\circ s)^{\ast}(TP))$ is a
vector bundle, we extend $u_{\phi}$ arbitrarily to the section $\widetilde
{u_{\phi}}:N\rightarrow\mathrm{Hom}(S^{r+1}((\pi_{N}^{\infty}\circ s)^{\ast
}(TN)),(\pi_{P}^{\infty}\circ s)^{\ast}(TP))$. Then we define $\widetilde{u}$
by $\widetilde{u}=\pi_{r}^{\infty}\circ s\oplus\widetilde{u_{\phi}}$ as the
section of $J^{r+1}(N,P)=$ $J^{r+1}(TN,TP)$. We lift $\widetilde{u}$ to the
section $s^{J}$\ of $J^{\infty}(N,P)$ over $N$. Then we have that $\pi
_{r+1}^{\infty}\circ s^{J}=\widetilde{u}$ and $\pi_{r}^{\infty}\circ s^{J}%
=\pi_{r}^{\infty}\circ s$.\ Furthermore, we define the homotopy $s_{\lambda
}\in\Gamma_{\Omega^{I}}(N,P)$ by%
\[
s_{\lambda}=(1-\lambda)s+\lambda s^{J}.
\]
It follows from $\pi_{r}^{\infty}\circ s_{\lambda}=\pi_{r}^{\infty}\circ
s=\pi_{r}^{\infty}\circ s^{J}$ that $s_{\lambda}$ is transverse to
$\Sigma^{I_{r}}(N,P)$ and $S^{I_{r}}(s_{\lambda})=S^{I_{r}}(s^{J})=S^{I_{r}%
}(s)$.

We prove that $s^{J}\in\Omega^{J}(N,P)$. For any point $c\in S^{I_{r}}(s)$,
let $U_{c}$ be a convex neighborhood of $c$ and let $k$\ and $y_{p}\ $be the
coordinates of $\exp_{N,c}((K_{r})_{c})$ and $\exp_{P,\pi_{P}^{\infty}\circ
s^{J}(c)}((\pi_{P}^{\infty}\circ s^{J})^{TP}(Q)_{c})$\ respectively. Let
$D_{k}$ denote the vector of the total tangent bundle $\mathbf{D}$ which
corresponds $k$ defined in [B, definition 1.6]. It follows from the definition
of $\mathbf{D}$ that%
\[%
\begin{array}
[c]{ll}%
(\bigcirc^{r+1}D_{k})y_{p}|_{s^{J}(c)}=\partial^{r+1}y_{p}/\partial
k^{r+1}(c)\neq0 & \text{for }c\in S^{I_{r}}(s).
\end{array}
\]
Then it follows from Lemma 7.2 (2) that
\[
\mathbf{d}_{r+1,s^{J}(c)}|\mathbf{K}_{r,s^{J}(c)}:\mathbf{K}_{r,s^{J}%
(c)}\rightarrow\mathbf{P}_{r,s^{J}(c)}\supset\text{Hom}(\bigcirc^{r}%
\mathbf{K}_{r,s^{J}(c)},\mathbf{Q}_{s(c)})
\]
is injective. Hence, we have that $s^{J}(S^{I_{r}}(s))\subset\Sigma^{J}(N,P)$.
Since $s^{J}(N\setminus(S^{I_{r}}(s))\subset\Omega^{I_{r-1},0}(N,P)$, the
assertion is proved. This proves the proposition.
\end{proof}

\begin{proof}
[Proof of Theorem 0.2]By the assumption, $j^{\infty}f$ is the section
$N\rightarrow\Omega^{I}(N,P)$. By Proposition 7.1, we have the section
$s^{J}:N\rightarrow\Omega^{J}(N,P)$ such that $\pi_{P}^{\infty}\circ s$ and
$\pi_{P}^{\infty}\circ s^{J}$ are homotopic. By Theorem 0.1 we obtain an
$\Omega^{J}$-regular map $g$ such that $j^{\infty}g$ and $s^{J}$ are
homotopic. This proves the assertion.
\end{proof}

\begin{corollary}
Let $n\geq p\geq2,$ and $N$ and $P$ be as above. Let $I=(n-p+1,1,1,1)$ and
$J=(n-p+1,1,1,0)$ such that $n-p$ and $r$ ($r\geq3$) are odd integers. Then if
$f:N\rightarrow P$ is an $\Omega^{I}$-regular map with $j^{\infty}f\in
\Gamma_{\Omega^{I}}^{tr}(N,P)$, then $f$ is homotopic to an $\Omega^{J}%
$-regular map $g:N\rightarrow P$ such that $j^{\infty}f$ and $j^{\infty}g$ are
homotopic in $\Gamma_{\Omega^{I}}^{tr}(N,P)$ and that $S^{I}(j^{\infty
}f)=S^{I}(j^{\infty}g)$.
\end{corollary}

This corollary proves the Chess conjecture ([C]). Sadykov[Sady] has actually
proved this corollary for $J=(n-p+1,1,0)$ in the case of $N$ and $P$ being orientable.

Let $\pi_{0}(X)$ be the arcwise connected components of $X$. Theorem 0.1
asserts that%
\[
(j_{\Omega^{I}})_{\ast}:\pi_{0}(C_{\Omega^{I}}^{\infty}(N,P))\rightarrow
\pi_{0}(\Gamma_{\Omega^{I}}(N,P))
\]
is surjective. However, $(j_{\Omega^{I}})_{\ast}$ is not necessarily
injective. Let $N=S^{2}$, $P=\mathbf{R}^{2}$ and $I=(1,0)$. Then we have by
[An3] that $\Omega^{1,0}(2,2)$ is homotopy equivalent to $SO(3)$. it follows
from [Ste, 36.4] that every two sections of $\Omega^{1,0}(S^{2},\mathbf{R}%
^{2})$ over $S^{2}$ are mutually homotopic. Namely, $\pi_{0}(\Gamma
_{\Omega^{I}}(N,P))$ consists of a single element. On the other hand, let
$f_{\lambda}:S^{2}\rightarrow\mathbf{R}^{2}$ be a homotopy of fold-maps.
Define $F:S^{2}\times\lbrack0,1]\rightarrow\mathbf{R}^{2}$ by $F(x,\lambda
)=f_{\lambda}(x)$ so that if $\lambda$ is sufficiently small, then
$F(x,\lambda)=f_{0}(x)$ and $F(x,1-\lambda)=f_{1}(x)$. By a very small
perturbation of $F$ fixing $f_{0}$ and $f_{1}$, we may assume that $F$ is
smooth and $f_{\lambda}$\ is still an $\Omega^{1,0}$-regular map for any
$\lambda$. Furthermore, the map $F:S^{2}\times\lbrack0,1]\rightarrow
\mathbf{R}^{2}\times\lbrack0,1]$ becomes an $\Omega^{1,0}$-regular map, and
$S^{1,0}(F)$ is a submanifold of $S^{2}\times\lbrack0,1]$. By the Jacobian
matrix of $F$ we know that the kernel line bundle $K_{1}(j^{\infty}F)$ over
$S^{1,0}(F)$ is independent with $\partial/\partial\lambda$, and
$T(S^{1,0}(F))\cap K_{1}(j^{\infty}F)=\{\mathbf{0}\}$. This implies that
$S^{1,0}(F)$ is regularly projected onto $[0,1]$. Hence, $S^{1,0}(f_{0})$ must
be diffeomorphic to $S^{1,0}(f_{1})$. Thus we conclude that $\pi_{0}%
(C_{\Omega^{I}}^{\infty}(S^{2},\mathbf{R}^{2}))$ is an infinite set.

\end{document}